\documentclass[10.5pt,a4paper]{article}

\usepackage[leqno]{amsmath}
\usepackage{graphicx,latexsym,euscript,makeidx,color,bm}
\usepackage{epstopdf,subfigure}
\usepackage{amsmath,amsfonts,amssymb,amsthm,thmtools,mathtools,mathrsfs,enumerate}
\usepackage[colorlinks,linkcolor=blue,citecolor=red]{hyperref}
\usepackage{geometry}
\allowdisplaybreaks[4]
\def\dbE{\mathbb{E}}    
\def\dbF{\mathbb{F}}  \def\cF{{\cal F}}  
    
\def\dbH{\mathbb{H}}

  \def\cL{{\cal L}}  
  \def\cM{{\cal M}}

\def\dbP{\mathbb{P}}    
    
\def\dbR{\mathbb{R}}

  \def\cY{{\cal Y}}  
  \def\cZ{{\cal Z}}  
\def\ss{\smallskip}             \def\hb{\hbox}
\def\ms{\medskip}              
\def\bs{\bigskip}            \def\as{\hb{a.s.}}
       \def\tr{\hb{tr$\,$}}
         
\def\no{\noindent}          
         
\def\nn{\nonumber}         
\def\rf{\eqref}            
\def\cd{\cdot}             
     \def\({\Big (}       \def\ba{\begin{aligned}}
\def\les{\leqslant}      \def\){\Big )}       \def\ea{\end{aligned}}
\def\ges{\geqslant}      \def\[{\Big[}        \def\bel{\begin{equation}\label}
\def\ti{\tilde}          \def\]{\Big]}        \def\ee{\end{equation}}
      \def\q{\quad}        
\def\h{\widehat}         \def\qq{\qquad}      
                                              
\def\a{\alpha}        \def\Om{\Omega}  
   \def\D{\Delta}   \def\d{\delta}        
         \def\si{\sigma}
\def\f{\varphi}   \def\l{\lambda}        \def\e{\varepsilon}
\def\t{\tau}    \def\i{\infty}      
\newtheoremstyle{thry}% name
{}      % Space above
{}      % Space below
{\sl}   % Body font
{}      % Indent amount
{\bf}   % Theorem head font
{.}     % Punctuation after theorem head
{.5em}  % Space after theorem head
{}      % Theorem head spec (can be left empty, meaning 'normal')
\theoremstyle{thry}

\newtheorem{theorem}{Theorem}[section]
\newtheorem{proposition}[theorem]{Proposition}

\newtheorem{lemma}[theorem]{Lemma}

\theoremstyle{definition}
\newtheorem{definition}[theorem]{Definition}

\theoremstyle{remark}
\newtheorem{remark}[theorem]{Remark}

\def\punct{}
\newtheoremstyle{dotless}{}{}{\rm}{}{\bf}{\punct}{.5em}{}
\theoremstyle{dotless}

\newtheorem{taggedassumptionx}{}
\newenvironment{taggedassumption}[1]
 {\renewcommand\thetaggedassumptionx{#1}\taggedassumptionx}
 {\endtaggedassumptionx}

\makeatletter
   
   \@addtoreset{equation}{section}
   \newcommand{\setword}[2]{%
   \phantomsection
   #1\def\@currentlabel{\unexpanded{#1}}\label{#2}%
   }
\makeatother

\begin{document}

\title{\bf  Feynman--Kac Representation for Parabolic
Equations with Memory: A New Class of BSDEs}

\author{
Hanxiao Wang\thanks{School of Mathematical Sciences, Shenzhen University, Shenzhen,
518060, China (Email: {\tt hxwang@szu.edu.cn}). HW is supported in part by   National Natural Science Foundation of China (Nos. 12522121 and 12631016).}
~~~~
Jiongmin Yong\thanks{Department of Mathematics, University of Central Florida, Orlando, FL 32816, USA
                 (Email: {\tt jiongmin.yong@ucf.edu}).
                  JY is supported in part by  Natural Science Foundation (No. DMS-2305475).}
~~~~
Qi Zhang\thanks{School of Mathematical Sciences and Laboratory of Mathematics for Nonlinear Science, Fudan University, Shanghai 200433, China (Email: {\tt qzh@fudan.edu.cn}). QZ is supported in part by National Key R\&D Program of China (No. 2022YFA1006101), and
National Natural Science Foundation of China (No. 12371445). }
}

\maketitle

%\linespread{1.5}

\no\bf Abstract.
\rm
This paper is concerned with the connection between parabolic equations with memory (MPEs, for short) and backward stochastic differential equations (BSDEs, for short). To provide a probabilistic representation of MPEs, we introduce a new class of BSDEs, called BSDEs with memory (MBSDEs, for short). The main feature of MBSDEs is that, at the current time $s$, the generator depends not only on the current value $Y^{t,x}(s)\equiv Y^{s,X^{t,x}(s)}(s)$ but also on the future path $r\mapsto Y^{r,X^{t,x}(s)}(r)$, $r\in[s,T]$. This structure reveals a fundamental difference from classical BSDEs: in the classical setting, the time and space variables are coupled and inseparable, whereas in MBSDEs they are completely independent. We prove the well-posedness of these MBSDEs and show that their adapted solutions provide a probabilistic representation for the solutions of MPEs. Furthermore, under a monotonicity condition, we establish a comparison principle for one-dimensional MBSDEs, which enables us to develop a viscosity solution theory for MPEs. This generalizes the celebrated nonlinear Feynman--Kac formula of Peng \cite{Peng1991} and Pardoux--Peng \cite{Pardoux-Peng1992} to the case with memory, in both the classical and viscosity solution senses.
\ms

\no\bf Keywords. \rm  parabolic equation with memory, backward stochastic differential equation,
 Feynman--Kac formula.
 %non-local operator, path-dependence,

\ms

\no\bf AMS subject classifications. \rm
45K05,  35K10, 60H10, 60H20.

\section{Introduction}
The nonlinear Feynman--Kac formula, established in the seminal works of Peng \cite{Peng1991} and Pardoux--Peng \cite{Pardoux-Peng1992}, asserts that the solution to the second-order semilinear parabolic equation
$$\left\{\begin{aligned}
&\partial_tu(t,x)+{1\over 2}\tr\big[\partial_{xx}^2 u(t,x)\si(t)\si(t)^\top\big]\\
&\q+g\big(t,x,u(t,x),\partial_x u(t,x)\si(t)\big)=0,\qq (t,x)\in[0,T]\times\dbR^n,\\
&u(T,x)=h(x),\qq x\in\dbR^n,
\end{aligned}\right.
$$
admits the following probabilistic representation:
$$
u(t,x)=Y^{t,x}(t),\qq (t,x)\in[0,T]\times\dbR^n,
$$
where $Y^{t,x}(\cd)$, together with some adapted process $Z^{t,x}(\cd)$, is the unique adapted solution to the following Markovian backward stochastic differential equation (BSDE, for short):
\bel{BSDE-Peng}\left\{\begin{aligned}
&dX^{t,x}(s)=\si(s) dW(s),\qq s\in[t,T],\\
&dY^{t,x}(s)=-g\big(s,X^{t,x}(s),Y^{t,x}(s),Z^{t,x}(s)\big)ds+Z^{t,x}(s)dW(s),\qq s\in[t,T],\\
&X^{t,x}(t)=x,\qq Y^{t,x}(T)=h(X^{t,x}(T)).\end{aligned}\right.
\ee
This connection extends the Feynman--Kac representation from linear parabolic PDEs to semilinear PDEs, thereby enabling a probabilistic treatment of such PDEs. Moreover, it paves the way for the use of Monte Carlo methods in their numerical solution; see e.g. Ma--Yong \cite{Ma-Yong-1999} and Delarue--Menozzi \cite{Delarue-Menozzi2006}. Since then, a substantial body of work has been devoted to studying the connections between forward--backward stochastic differential equations (FBSDEs, for short) and parabolic PDEs; see e.g. Ma--Protter--Yong \cite{Ma-Protter-Yong-1994}, Pardoux--Tang \cite{Pardoux-Tang1999}, and Delarue \cite{Delarue2002} for the quasi-linear case, and Cheridito--Soner--Touzi--Victoir \cite{Cheridito}, Peng \cite{Peng2010,Peng2019}, and Touzi \cite{Touzi2010} for the fully nonlinear case. Further developments on the connection between FBSDEs and stochastic PDE, path-dependent PDE, or mean-field PDEs can be found in Pardoux--Peng \cite{Pardoux-Peng1994}, Hu--Ma--Yong \cite{Hu-Ma-Yong2002}, Ekren--Christian--Touzi--Zhang \cite{Ekren--Christian--Touzi--Zhang2014}, and Buckdahn--Li--Peng--Rainer \cite{Buckdahn--Li--Peng--Rainer2017}, among others. Recently, this connection has played a crucial role in the study of the master equation in mean-field games and of deep-learning methods for high-dimensional PDEs; see Cardaliaguet--Delarue--Lasry--Lions \cite{Cardaliaguet--Delarue--Lasry--Lions2019} and Han--Jentzen--E \cite{Han2016}, respectively.

\ss
The parabolic equation with memory (MPE), a variant of the classical parabolic PDE, serves to model memory effects in dynamical systems. Typical examples include the equations arising in the electromagnetic theory of dielectrics and in heat conduction for materials with memory, which were investigated by Davis \cite{Davis1976,Davis1979} in the 1970s. For a detailed account of the background and motivations of MPEs, we refer to the book by Amendola--Fabrizio--Golden \cite{Amendola2021} and the references therein. The presence of the memory term makes the MPE fundamentally different from its classical counterpart. A representative example is that the flow generated by the heat equation with memory exhibits a hybrid parabolic--hyperbolic character; see Gurtin--Pipkin \cite{Gurtin} and Yong--Zhang \cite{Yong-Zhang2005,Yong-Zhang2011}. Due to this distinctive feature, the MPE has drawn considerable interest in recent decades. For instance, Alabau-Boussouira--Cannarsa--Sforza \cite{Alabau-2008} obtained a decay estimate; Fabrizio--Giorgi--Pata \cite{Fabrizio2010} proposed a new approach based on several novel notions; Andrade and Viana \cite{Andrade-2017} established well-posedness and regularity within an abstract framework; Wang--Zhang--Zuazua \cite{WZZ-2022} carried out a detailed analysis of the flow generated by MPEs, to name but a few. The associated control problems have also been extensively studied;
see e.g. Aassila--Cavalcanti--Soriano \cite{Aassila2000},  Fu--Yong--Zhang \cite{Fu2009}, Chaves-Silva--Zhang--Zuazua \cite{Chaves2017}, and Wang--Zhang--Zuazua \cite{WZZ-2024}.
%Some numerical results can be found in ......

\ss
Then a natural question arises concerning the \emph{probabilistic interpretation} of MPEs, namely the corresponding \emph{Feynman--Kac formula}. This question, however, remains unexplored, and the present paper aims to address this important issue. Several new and interesting phenomena emerge in the process. To be specific, we study the probabilistic representation of the following $m$-dimensional MPE:
\bel{H-PDE}\left\{\begin{aligned}
&\partial_tu(t,x)+{1\over 2}\tr\big[\partial_{xx}^2 u(t,x)\si(t)\si(t)^\top\big]\\
&\q+g\big(t,x,u(t,x),\partial_x u(t,x)\si(t),\cM[u](t,x)\big)=0,\qq (t,x)\in[0,T]\times\dbR^n,\\
&u(T,x)=h(x),\q x\in\dbR^n,
\end{aligned}\right.
\ee
where the memory operator $\cM[\,\cd\,]$ is defined by
$$
\cM[u](t,x)= \int_t^T M(s-t) u(s,x)ds,\qq (t,x)\in[0,T]\times\dbR^n,
$$
with $M(\cd):[0,T]\to\dbR^{k\times m}$ being a Lipschitz continuous function. Observe that, for any fixed point $(t,x)$, the quantity $\cM[u](t,x)$ depends on the values $u(s,x)$ for $s\in[t,T]$. Hence, $\cM[\,\cd\,]$ is a non-local operator in the time variable.

\ss
To address  this non-local term, we introduce a new class of BSDEs:
\bel{MBSDE}\left\{\begin{aligned}
&dX^{t,x}(s)=\si(s) dW(s),\qq s\in[t,T],\\
&dY^{t,x}(s)=-g\big(s,X^{t,x}(s),Y^{t,x}(s),Z^{t,x}(s),\cM^{t,x}[Y](s)\big)ds+Z^{t,x}(s)dW(s),\qq s\in[t,T],\\
&X^{t,x}(t)=x,\q Y^{t,x}(T)=h(X^{t,x}(T)),\end{aligned}\right.
\ee
where $(t,x)\in[0,T]\times\dbR^n$ is the initial pair, and for any fixed $s\in[t,T]$,
\bel{def-F}
\cM^{t,x}[Y](s)=\int_s^T M(r-s) Y^{r,X^{t,x}(s)}(r)dr,\qq s\in[t,T].
\ee
The key novelty of \rf{MBSDE} is the inclusion of the non-local term $\mathcal{M}^{t,x}[Y](\cdot)$ in the generator $g$ which is designed to capture the memory term $\mathcal{M}[u](\cdot,\cdot)$ of the MPE \rf{H-PDE}. The superscript of $Y^{t,x}(\cdot)$ signifies that the initial pair of the associated diffusion process is $(t,x)$. In the integral term \rf{def-F}, the variable $r$ plays a dual role: it is both the current time and the initial time of $Y^{r,X^{t,x}(s)}(r)$. One may therefore expect that $Y^{r,X^{t,x}(s)}(r)$ is $\mathcal{F}_s$-measurable for any $r\in[s,T]$, which ensures that the memory term \rf{def-F} remains progressively measurable. To distinguish this equation from classical BSDEs, we refer to it as a \emph{BSDE with memory} (\emph{MBSDE}, for short). The identification of this appropriate form constitutes the first contribution of our paper.

\ss
Recalling \rf{BSDE-Peng} and invoking the flow property of BSDEs, the solution to \rf{BSDE-Peng} satisfies
$$
Y^{t,x}(s)=Y^{s,X^{t,x}(s)}(s)=Y^{s,\ti x}(s)|_{\ti x=X^{t,x}(s)}.
$$
Consequently, in BSDE \rf{BSDE-Peng}, the space variable is always evaluated at the current value of the diffusion process $X^{t,x}(\cd)$ at time $s$, so the time and space variables are coupled and inseparable. The situation is entirely different for MBSDE \rf{MBSDE}: there we need to determine $Y^{r,X^{t,x}(s)}(r)=Y^{r,\ti x}(r)|_{\ti x=X^{t,x}(s)}$ for $s\neq r$. Hence, in MBSDEs the time and space variables are completely independent, marking an essential difference from the classical BSDE framework. This motivates the more natural notation $(Y^{\cd,\cd}(\cd),Z^{\cd,\cd}(\cd))$, rather than $(Y^{t,x}(\cd),Z^{t,x}(\cd))$, for the solution to \rf{MBSDE}.

\ss
We now turn to the well-posedness of the MBSDE \rf{MBSDE}. The main obstacle is the non-local term $\cM^{t,x}[Y](\cd)$. We use an approximation argument to deal with it and thereby construct a solution to \rf{MBSDE}. Specifically, we partition $[0,T]$ into $N$ subintervals $[t_i,t_{i+1}]$ with $t_i={i\over N}$, and on each subinterval we build a suitable standard BSDE to approximate the MBSDE. A key tool in this construction is the Feynman--Kac formula for backward stochastic Volterra integral equations (BSVIEs, for short) established in Wang--Yong--Zhang \cite{Wang-Yong-Zhang2022}. Gluing these local approximations together over $[0,T]$ yields an approximate solution to \rf{MBSDE} that depends on $N$. Through a series of suitable estimates, we prove that this sequence of approximate solutions is Cauchy in an appropriate norm, and its limit is exactly a solution to \eqref{MBSDE}. This approach not only establishes the existence of a solution to \eqref{MBSDE} but also reveals its connections with and differences from classical BSDEs.  Uniqueness then follows from a standard argument.

\ss
Given the unique solution $(Y^{\cd,\cd}(\cd),Z^{\cd,\cd}(\cd))$ to \eqref{MBSDE}, the following function is well-defined:
$$
\bar u(t,x)=Y^{t,x}(t),\qq (t,x)\in[0,T]\times\dbR^n.
$$
When the coefficients $\si(\cd)$, $g(\cd)$, and $h(\cd)$ are sufficiently smooth, we can show that $\bar u(\cd,\cd)$ belongs to $C^{1,2}([0,T]\times\dbR^n;\dbR^m)$, and hence is a classical solution to \eqref{MBSDE}. This establishes a rigorous connection between the MPE \rf{H-PDE} and the MBSDE \rf{MBSDE}. We emphasize that, unlike analytic methods, our approach does not require the leading coefficient $\si(\cd)$ to be nondegenerate. To handle the non-local term $\cM^{t,x}[Y](\cd)$, we need to establish the uniform estimate $\sup_{t,s,x}|Y^{t,x}(s)|$.

\ss
When the coefficients are less regular, we adopt the viscosity solution framework under the assumption that $Y^{\cd,\cd}(\cd)$ is scalar. Since viscosity solutions are local in nature whereas the MPE \rf{H-PDE} contains a non-local term, its comparison principle becomes considerably more delicate. Indeed, the comparison principle may fail in general, even for classical solutions. Inspired by Wang--Yong \cite{Wang-Yong2019} and Wang--Yong--Zhang \cite{Wang-Yong-Zhang2022}, we impose a monotonicity condition, under which the comparison principle holds naturally for classical solutions. The comparison principle for viscosity solutions is then established by a modification method.

\ss
The rest of the paper is organized as follows. Section \ref{Sec:Preliminaries} collects some preliminary results and introduces a few elementary notions. Section \ref{sec:3} is devoted to identifying the form of MBSDEs. In Section \ref{sec:4}, we establish the well-posedness and the comparison theorem for MBSDEs. The Feynman--Kac formula is proved in Section \ref{sec:6} for classical solutions and in Section \ref{sec:7} for viscosity solutions.

\section{Preliminaries}\label{Sec:Preliminaries}

Let $(\Om,\cF,\dbP)$ be a complete probability space, on which a standard $1$-dimensional Brownian motion
$W(\cd)$ is defined. The augmented natural filtration of $W(\cd)$ is denoted by $\dbF=\{\cF_t\}_{t\ges0}$.
Let $T>0$ be a fixed time horizon and $t\in[0,T)$. We denote
$$
\D[t,T]=\big\{(s,r)\in[t,T]^2\bigm|t\les s\les r\les T\big\}.
$$
We further introduce the following spaces of functions and processes:
\begin{align*}
L_\dbF^p(t,T;\dbR^m)
 &  =\Big\{\f:[t,T]\times\Om\to\dbR^m~\big|~ \f(\cd) \text{ is $\dbF$-progressively measurable,} \\
   &\qq\qq\qq \text{with } \dbE\Big[\Big(\int_t^T|\f(s)|^2ds\Big)^{p/2}\Big]<\i\Big\}; \\
L_\dbF^p(\D[t,T];\dbR^m)
 &  =\Big\{\f:\D[t,T]\times\Om\to\dbR^m~\big|~ \f(s,\cd)\in L_\dbF^p(s,T;\dbR^m)~ \as~ s\in[t,T], \\
   &\qq\qq\qq \text{with } \sup_{s\in[t,T]}\dbE\Big[\Big(\int_s^T|\f(s,r)|^2dr\Big)^{p/2}\Big]<\i\Big\}; \\
L_\dbF^p(\Om;C([t,T];\dbH))
&   =\Big\{\f:[t,T]\times\Om\to\dbR^m~\big|~ \f(\cd) \text{ is $\dbF$-adapted and continuous,} \\
   &\qq\qq\qq \text{with } \dbE\Big[\sup_{s\in[t,T]}|\f(s)|^p\Big]<\i\Big\};\\
C([0,T]\times\dbR^n;\dbR^m)
   &=\Big\{\f:[0,T]\times\dbR^n\to\dbR^m ~\big|~ \f(\cd) \text{ is a continuous function}\Big\}; \\
C^{i,j}([0,T]\times\dbR^n;\dbR^m)
  & =\Big\{\f(\cd,\cd)\in C([0,T]\times\dbR^n;\,\dbR^m) ~\big|~ \f(\cd,\cd) \text{ is $i$-times and $j$-times continuously } \\
   &\qq \text{differentiable in the time and space variables, respectively}\Big\}.
\end{align*}

\ss

Next, we recall the following results on BSVIEs which will play an important role in the study of MBSDEs. Consider the following BSVIE:
\begin{align}
\eta^{t,x}(s)&=l(s,X^{t,x}(T))+ \int_s^T f(s,r,X^{t,x}(r),\eta^{t,x}(r),\zeta^{t,x}(s,r))dr\nn\\
&\q-\int_s^T \zeta^{t,x}(s,r)dW(r),\qq s\in[t,T],\label{BSVIE}
\end{align}
where for any given initial pair $(t,x)\in[0,T)\times\dbR^n$, the process $X^{t,x}(\cd)$ is uniquely given by
\begin{align*}
X^{t,x}(s)=x+\int_t^s \si(r)dW(r),\qq s\in[t,T].\label{SDE}
\end{align*}

\ss
We impose the following assumption on  the BSVIE \rf{BSVIE}.

\begin{taggedassumption}{(H0)}\label{ass:H0}
The mappings $f:\D[0,T]\times\dbR^n\times\dbR^m\times\dbR^{m}\to\dbR^m$, $l:[0,T]\times\dbR^n\to\dbR^m$, and $\si:[0,T]\to\dbR^n$ are jointly measurable. There exists a constant $L>0$ such that
\begin{align*}
&|\si(s)|+|l(0,0)|+|f(0,s,0,0,0)|\les L,\qq \forall s\in [0,T];\\
&|l(s_1,x_1)-l(s_2, x_2)|+|f(s_1,r,x_1,\eta_1,\zeta_1)-f(s_2,r,x_2,\eta_2,\zeta_2)|\\
&\q\les L(|s_1-s_2|+|x_1-x_2|+|\eta_1-\eta_2|+|\zeta_1-\zeta_2|),\\
&\qq~\forall (s_i,r,x_i,\eta_i,\zeta_i)\in [0,T]^2\times\dbR^n\times\dbR^m\times\dbR^m,\q i=1,2.
\end{align*}
\end{taggedassumption}

From Wang--Yong--Zhang \cite{Wang-Yong-Zhang2022}, we have the following result.

\begin{proposition}\label{prop:bsvie}
Under {\rm\ref{ass:H0}},  for any $p>1$ and initial pair $(t,x)\in[0,T)\times\dbR^n$, the BSVIE \rf{BSVIE} admits a unique solution
$(\eta^{t,x}(\cd),\zeta^{t,x}(\cd,\cd))\in L^p_{\dbF}(\Om; C([t,T];\dbR^m))\times  L^p_{\dbF}([t,T];\dbR^m)$,
satisfying the estimate
$$
\dbE_t\[\sup_{s\in[t,T]}|\eta^{t,x}(s)|^p\]+\sup_{s\in[t,T]}\dbE_t\[\int_s^T |\zeta^{t,x}(s,r)|^2dr\]^{p\over 2}
\les K_p\dbE_t\[1+\sup_{s\in[t,T]}|X^{t,x}(s)|^p\],
$$
for some constant $K_p>0$.
Moreover, the following deterministic and jointly continuous function is well-defined:
$$
u^0(t,x)=\eta^{t,x}(t), \qq (t,x)\in[0,T]\times\dbR^n,
$$
and it is  Lipschitz continuous in the state variable $x$.
\end{proposition}

\section{From MPEs to MBSDEs}\label{sec:3}
Recall that the MPE \rf{H-PDE} reads
\bel{H-PDE-rw}\left\{\begin{aligned}
&\partial_tu(t,x)+{1\over 2}\tr\big[\partial_{xx}^2 u(t,x)\si(t)\si(t)^\top\big]\\
&\q+g\big(t,x,u(t,x),\partial_x u(t,x)\si(t),\cM[u](t,x)\big)=0,\qq (t,x)\in[0,T]\times\dbR^n,\\
&u(T,x)=h(x),\qq x\in\dbR^n,
\end{aligned}\right.
\ee
where the memory operator $\cM[\,\cd\,]$ is defined by
\bel{def-M-operator-rw}
\cM[u](t,x)= \int_t^T M(s-t) u(s,x)ds,\qq (t,x)\in[0,T]\times\dbR^n.
\ee

\ss
We next derive the associated stochastic equation.
Since the leading coefficient of \rf{H-PDE-rw} is given by ${1\over2}\si\si^\top$,
the corresponding diffusion process naturally takes the form:
$$\left\{\begin{aligned}
&dX^{t,x}(s)=\si(s) dW(s),\qq s\in[t,T],\\
&X^{t,x}(t)=x.\end{aligned}\right.
$$
Assume the MPE \rf{H-PDE-rw}--\rf{def-M-operator-rw} admits a classical solution $u(\cd,\cd)$.
Applying It\^{o}'s formula to the mapping $s\mapsto u(s,X^{t,x}(s))$, we obtain
\begin{align}
d u(s,X^{t,x}(s))&=\(\partial_s u(s,X^{t,x}(s))+{1\over 2}\tr\big[\partial_{xx}^2 u(s,X^{t,x}(s))\si(s)\si(s)^\top\big]\)ds\nn\\
&\q +\partial_x u(s,X^{t,x}(s))\si(s)dW(s)\nn\\
&=-g\big(s,X^{t,x}(s),u(s,X^{t,x}(s)),\partial_x u(s,X^{t,x}(s))\si(s),\cM[u](s,X^{t,x}(s))\big)ds\nn\\
&\q+\partial_x u(s,X^{t,x}(s))\si(s)dW(s),\qq s\in[t,T],\label{proof-PDE-to-BSDE1}
\end{align}
with
\bel{proof-PDE-to-BSDE2}
\cM[u](s,X^{t,x}(s))=\int_s^T M(r-s)u(r,X^{t,x}(s))dr.
\ee
Motivated by Peng \cite{Peng1991}, it is natural to define the following two adapted processes:
\bel{proof-PDE-to-BSDE3}
Y^{t,x}(s)=u(s,X^{t,x}(s)),\qq Z^{t,x}(s)=\partial_xu(s,X^{t,x}(s))\si(s),\qq s\in[t,T].
\ee
Then we can rewrite \rf{proof-PDE-to-BSDE1} as follows:
\begin{align}
dY^{t,x}(s)
&=-g\big(t,x,Y^{t,x}(s),Z^{t,x}(s),\cM[u](s,X^{t,x}(s))\big)ds+Z^{t,x}(s)dW(s),\qq s\in[t,T].\label{proof-PDE-to-BSDE4}
\end{align}
Given the above, the key is to represent $\cM[u](\cd,X^{t,x}(\cd))$ as an adapted  process.
From \rf{proof-PDE-to-BSDE3}, we observe that
$$
u(t,x)=Y^{t,x}(t),\qq  (t,x)\in[0,T]\times\dbR^n.
$$
Thus, the following relationship is expected to hold:
$$
u(r,X^{t,x}(s))=Y^{r,\ti x}(r)|_{\ti x=X^{t,x}(s)}=: Y^{r,X^{t,x}(s)}(r)  ,\qq t\les s\les r\les T,
$$
which, together with \rf{proof-PDE-to-BSDE2}, leads to
$$
\cM[u](s,X^{t,x}(s))=\int_s^T M(r-s) Y^{r,X^{t,x}(s)}(r)dr=:\cM^{t,x}[Y](s),\qq s\in[t,T].
$$
Substituting the above into \rf{proof-PDE-to-BSDE4},
the stochastic equation associated with MPE \rf{H-PDE-rw}--\rf{def-M-operator-rw} is derived as follows:
\bel{MBSDE-rw}\left\{\begin{aligned}
&dX^{t,x}(s)=\si(s) dW(s),\qq s\in[t,T],\\
&dY^{t,x}(s)=-g\big(s,X^{t,x}(s),Y^{t,x}(s),Z^{t,x}(s),\cM^{t,x}[Y](s)\big)ds+Z^{t,x}(s)dW(s),\qq s\in[t,T],\\
&X^{t,x}(t)=x,\q Y^{t,x}(T)=h(X^{t,x}(T)),\end{aligned}\right.
\ee
where  $(t,x)\in[0,T]\times\dbR^n$ is the initial pair, and for any fixed $s\in[t,T]$,
\bel{def-F-rw}
\cM^{t,x}[Y](s)=\int_s^T M(r-s) Y^{r,X^{t,x}(s)}(r)dr,\qq s\in[t,T].
\ee

\ss
\begin{remark}
The superscript of $Y^{t,x}(\cd)$ indicates that the initial pair of the corresponding diffusion process is $(t,x)$. Note that for $Y^{r,X^{t,x}(s)}(r)$, the initial state of the diffusion process at time $r$ is given by $X^{t,x}(r)$, which is both $\cF_r$-measurable and $\cF_s$-measurable. Consequently, $Y^{r,X^{t,x}(s)}(r)$ is expected to be $\cF_s$-measurable for any $r\in[s,T]$, which ensures that the memory term \rf{def-F-rw} remains progressively measurable.
\end{remark}

\begin{remark}
%Recalling \rf{BSDE-Peng} and using the flow property of FBSDEs, the solution to \rf{BSDE-Peng} satisfies
%$$
%Y^{t,x}(s)=Y^{s,X^{t,x}(s)}(s)=Y^{s,\ti x}(s)|_{\ti x=X^{t,x}(s)}.
%$$
%Thus, in BSDE \rf{BSDE-Peng}, the space variable is always taken as the current value $X^{t,x}(s)$ of the forward process, so the time and space variables are coupled and inseparable. In %contrast, for the equation \rf{MBSDE-rw}, we need to determine the value
%$$
%Y^{r,X^{t,x}(s)}(r)=Y^{r,\ti x}(r)|_{\ti x=X^{t,x}(s)}
%$$
%with $s\neq r$. Hence,
Notice that the time and space variables of $Y^{r,X^{t,x}(s)}(r)$ in \rf{MBSDE-rw} are completely independent since $r$ is not necessarily equal to $s$, which is essentially different from the classical BSDE. In other words, in BSDE \rf{BSDE-Peng}, $(t,x)$ serves merely as a parameter, whereas in \rf{MBSDE-rw}, it truly functions as a variable. It is therefore more natural to denote the solution to \rf{MBSDE-rw} by $(Y^{\cd,\cd}(\cd),Z^{\cd,\cd}(\cd))$ instead of $(Y^{t,x}(\cd),Z^{t,x}(\cd))$.
%To distinguish this type of equation  from classical BSDEs, we call it a \emph{BSDE with memory} (\emph{MBSDE}, for short).
\end{remark}

\begin{definition}
We call $(Y^{\cd,\cd}(\cd),Z^{\cd,\cd}(\cd))$  an adapted solution to the MBSDE \rf{MBSDE-rw} if
\begin{enumerate}[(i)]
\item for any given $(t,x)\in[0,T]\times\dbR^n$, $ (Y^{t,x}(\cd),Z^{t,x}(\cd))\in L^2_{\dbF}(\Om; C([t,T];\dbR^m))\times L^2_{\dbF}([t,T];\dbR^m) $;
\item for any given  $0\les t\les s\les r\les T$ and $x\in\dbR^n$,
$Y^{r,X^{t,x}(s)}(r)= Y^{r,\ti x}(r)|_{\ti x=X^{t,x}(s)}$, and it is $\cF_s$-measurable;
\item    \rf{MBSDE-rw} is satisfied in the usual It\^{o}'s sense for any given $(t,x)\in[0,T]\times\dbR^n$.
\end{enumerate}
\end{definition}

From the above procedure of deriving the MBSDE from MPE, we obtain the following result.

\begin{theorem}\label{thm:PDE-to=BSDE}
Suppose that  the MPE \rf{H-PDE-rw}--\rf{def-M-operator-rw} admits a classical solution $u(\cd,\cd)$.
Then the process $(Y^{\cd,\cd}(\cd),Z^{\cd,\cd}(\cd))$,  defined by
\begin{align*}
\big(Y^{t,x}(s),\,Z^{t,x}(s)\big):=\big(u(s,X^{t,x}(s)),\,u_x(s,X^{t,x}(s))\sigma(s)\big),\qq
 s\in[t,T],~(t,x)\in[0,T]\times\dbR^n,
\end{align*}
is an adapted solution to  the MBSDE \rf{MBSDE-rw}.
\end{theorem}

The above shows that, given a classical solution to the MPE \rf{H-PDE-rw}--\rf{def-M-operator-rw} (if it exists), one can construct a solution to the MBSDE \rf{MBSDE-rw}. In the literature \cite{Peng1991,Ma-Zhang}, this is referred to as a \emph{representation theorem}. Conversely, the use of the MBSDE \rf{MBSDE-rw} solution to represent, in some sense, the solution to the MPE \rf{H-PDE-rw}--\rf{def-M-operator-rw} is more important and interesting, known as the \emph{Feynman--Kac formula}. We will establish this result in the rest of the paper.

\section{Well-posedness and comparison theorem of MBSDEs}\label{sec:4}

To establish the Feynman--Kac formula for the MPE \rf{H-PDE-rw}--\rf{def-M-operator-rw}, we first prove the well-posedness of the MBSDE \rf{MBSDE-rw}
under the following assumption.

\begin{taggedassumption}{(H1)}\label{ass:H1}
The mappings $g:[0,T]\times\dbR^n\times\dbR^m\times\dbR^{m}\times\dbR^k\to\dbR^m$, $h:\dbR^n\to\dbR^m$,
 $M:[0,T]\to\dbR^{k\times m}$,  and $\si:[0,T]\to\dbR^n$ are  jointly measurable.
 There exists a constant $L>0$ such that
\begin{align*}
&|\si(t)|+|h(x)|+|g(t,x,0,0,0)|\les L,\qq \forall (t,x)\in [0,T]\times\dbR^n;\\
&|h(x_1)-h(x_2)|+|g(s,x_1,y_1,z_1,y^\prime_1)-g(s,x_2,y_2,z_2,y^\prime_2)|+|M(t_1)-M(t_2)|\\
&\q\les L(|t_1-t_2|+|x_1-x_2|+|y_1-y_2|+|z_1-z_2|+|y^\prime_1-y^\prime_2|),\\
&\qq\forall (t_i,s,x_i,y_i,z_i,y_i^\prime)\in [0,T]^2\times\dbR^n\times\dbR^m\times\dbR^m\times\dbR^k,\,i=1,2.
\end{align*}
\end{taggedassumption}

\begin{theorem}\label{thm:well-posdeness}
Under {\rm\ref{ass:H1}}, the MBSDE   \rf{MBSDE-rw} admits a unique adapted solution $(Y^{\cd,\cd}(\cd),Z^{\cd,\cd}(\cd)) $ with the estimate:
\begin{align*}
&\sup_{t,s,x} \Big\{|Y^{t,x}(s)|^p +\dbE_s\[\int_s^T |Z^{t,x}(r)|^2dr\]^{p\over 2}\Big\}\nn\\
&\q\les K_p \sup_{t,s,x}  \dbE_s\[|h(X^{t,x}(T))|^p+\(\int_s^T |g(r,X^{t,x}(r),0,0,0)|dr\)^p \],
\end{align*}
where $K_p>0$ is a constant depending on $p>1$.
\end{theorem}

The proof of Theorem \ref{thm:well-posdeness} is presented below using an approximation approach.

\subsection{Approximation of MBSDEs}\label{subsec:approximate}

Since the appearance of $Y^{r,X^{t,x}(s)}(r)$, the MBSDE  \rf{MBSDE-rw}
has a special structure. We will prove its solvability  by an approximation approach.
A foundation is the Feynman--Kac formula. % (see \cite{Wang-Yong-Zhang2022}) for backward
%stochastic Volterra integral equations (BSVIEs, for short).

\ss

Note that the diffusion process given in \rf{MBSDE-rw} is
$$\left\{\begin{aligned}
&dX^{t,x}(s)=\si(s)dW(s),\qq s\in[t,T],\\
&X^{t,x}(t)=x,\end{aligned}\right.
$$
with $(t,x)\in[0,T)\times\dbR^n$.
We now give an approximating equation to the backward equation in \rf{MBSDE-rw}.
Let $\Pi_N\equiv\{t_k~|~k=0,1,..., N\}$ be a partition of $[0,T]$ with
$$
t_k={k\over N}T.
$$

\ss
\textbf{\emph{Step 1.}}
On $[t_{N-1},t_N]$, consider the following equation:
\bel{BSVIE-N}
\left\{\begin{aligned}
&Y_{N}^{t,x}(s)=h(X^{t,x}(T))+\int_s^T g\big(r,X^{t,x}(r),Y_N^{t,x}(r),Z_N^{t,x}(r),\cY_N^{t,x}(r)\big)dr-\int_s^T Z_N^{t,x}(r)dW(r),\\
&\cY_N^{t,x}(s)=\int_s^T M(r-s)Y_N^{t,x}(r)dr-\int_s^T\cZ_N^{t,x}(s,r)dW(r).
\end{aligned}\right.
\ee
%
%The first one in \rf{BSVIE-N} is a BSDE and the second one is a BSVIE.
As a whole, \rf{BSVIE-N} is a BSVIE. Under \ref{ass:H1}, by Proposition \ref{prop:bsvie},
the BSVIE \rf{BSVIE-N} admits a unique adapted solution $(Y_N^{t,x}(\cd),\cY_N^{t,x}(\cd),Z_N^{t,x}(\cd,\cd),
\cZ_N^{t,x}(\cd,\cd))$.
Moreover,
the following deterministic function is well-defined:
$$
V_N(s,x):=Y^{s,x}_N(s),\qq (s,x)\in[t_{N-1},T]\times\dbR^n,
$$
with $V_N(\cd,\cd)\in C^{0,1}$. Then the following is well-defined
\bel{YN-def}
Y^{r,X^{t,x}(s)}_N(r):=V_N(r,X^{t,x}(s))=Y^{s,\ti x}_N(s)|_{\ti x=X^{t,x}(s)},\qq (s,r,x)\in[t,T]\times[t_{N-1}\vee s,T]\times\dbR^n.
\ee
We remark that $Y^{r,X^{t,x}(s)}_N(r)$ is $\cF_s$-measurable with $s\les r$.

\ss
\textbf{\emph{Step 2.}}
We next consider the approximating equation over $[t_{N-2},t_N]$. Consider the following equation:
\bel{BSVIE-(N-1)'}\left\{\begin{aligned}
Y_{N-1}^{t,x}(s)&=h(X^{t,x}(T))+\int_s^T g\big(r,X^{t,x}(r),Y_{N-1}^{t,x}(r),Z_{N-1}^{t,x}(r),\cY_{N-1}^{t,x}(r)\big)dr\\
&\q-\int_s^TZ_{N-1}^{t,x}(r)dW(r),\qq s\in[t_{N-2},T],\\
\cY_{N-1}^{t,x}(s)&=\int_{s\vee t_{N-1}}^T M(r-s)Y^{r,X^{t,x}(s)}_N(r)dr-\int_s^T\cZ_{N-1}^{t,x}(s,r)dW(r),\qq s\in[t_{N-2},T].\end{aligned}\right.
\ee
Note that $Y^{r,X^{t,x}(s)}_N(r)$ has been well-defined by \rf{YN-def}, so the above is a (classical) BSVIE,
which admits a unique solution $(Y_{N-1}^{t,x}(\cd),\cY_{N-1}^{t,x}(\cd),Z_{N-1}^{t,x}(\cd),\cZ_{N-1}^{t,x}(\cd,\cd))$.
Moreover, using the fact that $Y^{r,X^{t,x}(s)}_N(r)\in\cF_s$ for any $r\in[s\vee t_{N-1},T]$, which implies that
$$
\int_{s\vee t_{N-1}}^T M(r-s)Y^{r,X^{t,x}(s)}_N(r)dr\in\cF_s,\qq s\in[t_{N-2},T],
$$
we have
$$
\cZ_{N-1}^{t,x}(\cd,\cd)\equiv 0.
$$
Then we can rewrite \rf{BSVIE-(N-1)'} as follows:
$$\left\{\begin{aligned}
Y_{N-1}^{t,x}(s)&=h(X^{t,x}(T))+\int_s^T g\big(r,X^{t,x}(r),Y_{N-1}^{t,x}(r),Z_{N-1}^{t,x}(r),\cY_{N-1}^{t,x}(r)\big)dr\\
&\q-\int_s^TZ_{N-1}^{t,x}(r)dW(r),\qq s\in[t_{N-2},T],\\
\cY_{N-1}^{t,x}(s)&=\int_{s\vee t_{N-1}}^T M(r-s)Y^{r,X^{t,x}(s)}_N(r)dr,\qq s\in[t_{N-2},T].\end{aligned}\right.
$$
By Proposition \ref{prop:bsvie} again,
the following  function is well-defined
$$
V_{N-1}(s,x):=Y^{s,x}_{N-1}(s),\qq (s,x)\in[t_{N-2},T]\times\dbR^n,
$$
with $V_{N-1}(\cd,\cd)\in C^{0,1}$.
Define
\bel{V-N-N-1}
V_{\Pi_N}(s,x)
:=\left\{\begin{aligned}
&V_N(s,x),\qq s\in[t_{N-1},T],\q x\in\dbR^n;\\
&V_{N-1}(s,x),\qq s\in[t_{N-2},t_{N-1}),\q x\in\dbR^n,\end{aligned}\right.
\ee
and
$$
Y^{t,x}_{\Pi_N}(s):=\left\{\begin{aligned}
&Y^{t,x}_N(s),\qq s\in[t_{N-1},T];\\
&Y^{t,x}_{N-1}(s),\qq s\in[t_{N-2},t_{N-1}),
\end{aligned}\right.
$$
$$
\cY^{t,x}_{\Pi_N}(s):=\left\{\begin{aligned}
&\cY^{t,x}_N(s),\qq s\in[t_{N-1},T];\\
&\cY^{t,x}_{N-1}(s),\qq s\in[t_{N-2},t_{N-1}).
\end{aligned}\right.
$$
Then the following relationship holds:
\begin{align*}
Y^{r,X^{t,x}(s)}_{\Pi_N}(r)=V_{\Pi_N}(r,X^{t,x}(s))=Y^{r,\ti x}_{\Pi_N}(r)|_{\ti x=X^{t,x}(s)},\q\q
 (s,r,x)\in[t,T]\times[t_{N-2}\vee s,T]\times\dbR^n.
\end{align*}

\textbf{\emph{Step 3.}}
We next consider the approximating equation over $[t_{N-3},t_N]$. Consider the following equation:
$$\left\{\begin{aligned}
Y_{N-2}^{t,x}(s)&=h(X^{t,x}(T))+\int_s^T g\big(r,X^{t,x}(r),Y_{N-2}^{t,x}(r),Z_{N-2}^{t,x}(r),\cY_{N-2}^{t,x}(r)\big)dr\\
&\q-\int_s^TZ_{N-2}^{t,x}(r)dW(r),\qq s\in[t_{N-3},T],\\
\cY_{N-2}^{t,x}(s)&=\int_{s\vee t_{N-2}}^T
M(r-s)Y^{r,X^{t,x}(s)}_{\Pi_N}(r)dr,\qq s\in[t_{N-3},T],
\end{aligned}\right.
$$
which admits a unique solution $(Y_{N-2}^{t,x}(\cd),\cY_{N-2}^{t,x}(\cd),Z_{N-2}^{t,x}(\cd))$.
Then the following  function is well-defined:
$$
V_{N-2}(s,x):=Y^{s,x}_{N-2}(s),\qq (s,x)\in[t_{N-3},T]\times\dbR^n.
$$
With the above, we can extend $V_{\Pi_N}(\cd,\cd)$ from $[t_{N-2},T]\times\dbR^n$ to $[t_{N-3},T]\times\dbR^n$ by setting
%
%\bel{V-N-N-2}V_{\Pi_N}(s,x):=\left\{\begin{aligned}
%%
%&V_{\Pi_N}(s,x),\qq s\in[t_{N-2},T],\q x\in\dbR^n;\\
%%
%&V_{N-2}(s,x),\qq s\in[t_{N-3},t_{N-2}),\q x\in\dbR^n,
%\end{aligned}\right.
$$V_{\Pi_N}(s,x):=
V_{N-2}(s,x),\qq s\in[t_{N-3},t_{N-2}),\q x\in\dbR^n,
$$
in which  $V_{\Pi_N}(s,x),\, (s,x)\in[t_{N-2},T]\times\dbR^n$ has been given by \rf{V-N-N-1}.
Similarly, we can extend  $Y^{t,x}_{\Pi_N}(\cd)$ from $[t_{N-2},T]\times\dbR^n$ to $[t_{N-3},T]\times\dbR^n$ by setting
%
%\bel{Y-N-N-2}
%Y^{t,x}_{\Pi_N}(s):=\left\{\begin{aligned}
%%
%&Y^{t,x}_{\Pi_N}(s),\qq s\in[t_{N-2},T];\\
%%
%%
%&Y^{t,x}_{N-2}(s),\qq s\in[t_{N-3},t_{N-2}),
%\end{aligned}\right.
$$
Y^{t,x}_{\Pi_N}(s):=Y^{t,x}_{N-2}(s),\qq s\in[t_{N-3},t_{N-2}),
$$
and extend  $\cY^{t,x}_{\Pi_N}(\cd)$ from $[t_{N-2},T]\times\dbR^n$ to $[t_{N-3},T]\times\dbR^n$ by setting
%
%\bel{cY-N-N-2}
%\cY^{t,x}_{\Pi_N}(s):=\left\{\begin{aligned}
%%
%&\cY^{t,x}_{\Pi_N}(s),\qq s\in[t_{N-2},T];\\
%%
%%
%&\cY^{t,x}_{N-2}(s),\qq s\in[t_{N-3},t_{N-2}),
%\end{aligned}\right.
%\ee
$$
\cY^{t,x}_{\Pi_N}(s):=\cY^{t,x}_{N-2}(s),\qq s\in[t_{N-3},t_{N-2}),
$$
with
\begin{align*}
Y^{r,X^{t,x}(s)}_{\Pi_N}(r)=V_{\Pi_N}(r,X^{t,x}(s))=Y^{r,y}_{\Pi_N}(r)|_{y=X^{t,x}(s)},\qq
 (s,r,x)\in[t,T]\times[t_{N-3}\vee s,T]\times\dbR^n.
\end{align*}

\ss
\textbf{\emph{Step 4.}}
Finally, we consider the general case. For each $k=0,1,2,...,N-1$, suppose that
$$
Y^{r,X^{t,x}(s)}_{\Pi_N}(r)=V_{\Pi_N}(r,X^{t,x}(s)),\qq (s,r,x)\in[t,T]\times[t_k\vee s,T]\times\dbR^n,
$$
has been well-defined.
Then on $[t_{k-1},T]$, we consider the following BSVIE:
\bel{BSVIE-k}
\left\{\begin{aligned}
Y_{k}^{t,x}(s)&=h(X^{t,x}(T))+\int_s^T g\big(r,X^{t,x}(r),Y_k^{t,x}(r),Z_k^{t,x}(r),\cY_k^{t,x}(r)\big)dr-\int_s^TZ_{k}^{t,x}(r)dr,\\
\cY_{k}^{t,x}(s)&=\int_{s\vee t_{k}}^T M(r-s) Y^{r,X^{t,x}(s)}_{\Pi_N}(r)dr.
\end{aligned}\right.
\ee
With the unique adapted solution $(Y_{k}^{t,x}(\cd),Y_{k}^{t,x}(\cd),\cY_{k}^{t,x}(\cd))$, we can define
$$
V_{k}(s,x):=Y^{s,x}_{k}(s),\qq (s,x)\in[t_{k},T]\times\dbR^n,
$$
and then extend $V_{\Pi_N}(\cd,\cd)$, $Y^{t,x}_{\Pi_N}(\cd)$ and $\cY^{t,x}_{\Pi_N}(\cd)$ from $[t_{k},T]\times\dbR^n$ to $[t_{k-1},t_{k}]\times\dbR^n$ by setting
%
%$$
%V_{\Pi_N}(s,x):=\left\{\begin{aligned}
%%
%&V_{\Pi_N}(s,x),\qq s\in[t_{k},T],\q x\in\dbR^n;\\
%%
%&V_{k}(s,x),\qq s\in[t_{k-1},t_{k}),\q x\in\dbR^n,
%\end{aligned}\right.
%$$
%%
%\bel{Y-N-k}
%Y^{t,x}_{\Pi_N}(s):=\left\{\begin{aligned}
%%
%&Y^{t,x}_{\Pi_N}(s),\qq s\in[t_k,T];\\
%%
%%
%&Y^{t,x}_{k}(s),\qq s\in[t_{k-1},t_k),
%\end{aligned}\right.
%\ee
%%
%%
%$$
%\cY^{t,x}_{\Pi_N}(s):=\left\{\begin{aligned}
%%
%&\cY^{t,x}_{\Pi_N}(s),\qq s\in[t_k,T];\\
%%
%%
%&\cY^{t,x}_k (s),\qq s\in[t_{k-1},t_k),
%\end{aligned}\right.
%$$
\begin{align}\label{Y-N-k}
&V_{\Pi_N}(s,x):=V_{k}(s,x),\qq s\in[t_{k-1},t_{k}),\q x\in\dbR^n,\nn\\
&Y^{t,x}_{\Pi_N}(s):=Y^{t,x}_{k}(s),\qq s\in[t_{k-1},t_k),\\
&\cY^{t,x}_{\Pi_N}(s):=\cY^{t,x}_k (s),\qq s\in[t_{k-1},t_k),\nn
\end{align}
with the relationship
\begin{align}
Y^{r,X^{t,x}(s)}_{\Pi_N}(r)=V_{\Pi_N}(r,X^{t,x}(s))=Y^{r,y}_{\Pi_N}(r)|_{y=X^{t,x}(s)},\ \ \
 (s,r,x)\in[t,T]\times[t_{k-1}\vee s,T]\times\dbR^n.\nn
\end{align}
We complete the construction of the solution when $k=0$.

\subsection{Convergence of approximating sequence }

%We first show that  the approximating sequence  $\{Y_{\Pi_N}^{t,x}(t);\, (t,x)\in[0,T]\times\dbR^n\}$ constructed by \rf{Y-N-k} has a limit as $N\to \infty$.

\begin{proposition}
Under {\rm\ref{ass:H1}}, the sequence  $\{Y_{\Pi_N}^{t,x}(t)\,|\, (t,x)\in[0,T]\times\dbR^n\}_{N=1}^\infty$ constructed by \rf{Y-N-k}
is Cauchy in the uniform norm $\|\cd\|_{C([0,T]\times\dbR^n;\dbR^m)}$. %, where $\Pi_N: t_0=0,t_1={T\over N},......,t_{N-1}={(N-1)T\over N},t_N=T$.
\end{proposition}

\begin{proof}
We divided the proof into three steps.

\ss
\emph{\textbf{Step 1.} A uniform  estimate for $\sup_{t,s\in\D[0,T],x\in\dbR^n}|Y_{\Pi_N}^{t,x}(s)|^2$.}

\ss

For any $(t,s)\in\Delta[0,T]$ and $x\in\dbR^n$ with $s\in[t_{k-1},t_k)$ for some $k$,
by the standard $L^2$-estimates for the BSDE \rf{BSVIE-k}, we have
\begin{align*}
\dbE_s\[\sup_{r\in[s,T]}|Y_{k}^{t,x}(r)|^2\]
&\les K \dbE_s\[|h(X^{t,x}(T))|^2+\int_s^T\int^T_{r\vee t_k}|Y_{\Pi_N}^{\t,X^{t,x}(r)}(\t)|^2 d\tau dr\\
&\qq\qq +\int_s^T |g(r,X^{t,x}(r),0,0,0)|^2dr \].
\end{align*}
Here and in the rest of this paper $K$ is a generic constant independent of $N$. Note that
$$
|Y_{\Pi_N}^{\t,X^{t,x}(r)}(\t)|^2\les\sup_{x\in\dbR^n} |Y_{\Pi_N}^{\t,x}(\t)|^2,
$$
and $\sup_{x\in\dbR^n} |Y_{\Pi_N}^{\t,x}(\t)|^2$ is a deterministic function.
Taking $r=s$ in the above, we get
\begin{align*}
|Y_{k}^{t,x}(s)|^2 &\les  K \dbE_s\[|h(X^{t,x}(T))|^2
+\int_s^T |g(r,X^{t,x}(r),0,0,0)|^2dr +\int_s^T \sup_{x\in\dbR^n}|Y_{\Pi_N}^{r,x}(r)|^2dr \]\\
&\les  K+K\int_s^T \sup_{x\in\dbR^n}|Y_{\Pi_N}^{r,x}(r)|^2dr.
\end{align*}
Noting the definition \rf{Y-N-k} for $Y_{\Pi_N}^{t,x}(s)$, we have
\begin{align*}
|Y_{\Pi_N}^{t,x}(s)|^2\les
K+K\int_s^T \sup_{x\in\dbR^n}|Y_{\Pi_N}^{r,x}(r)|^2dr,\qq s\in[0,T].
\end{align*}
Then by the Gr\"{o}nwall's inequality, we get the following uniform estimate:
\begin{align*}
\sup_{t,s\in\D[0,T],x\in\dbR^n}|Y_{\Pi_N}^{t,x}(s)|^2\les K.
\end{align*}
%in which $K>0$ is independent of $\Pi_N$.

\ss

\emph{\textbf{Step 2.} A uniform  estimate for $\sup_{t,s,x}|\nabla Y_{\Pi_N}^{t,x}(s)|^2$.}

\ss
Note that $ Y_{\Pi_N}^{t,x}(\cd)$ is constructed by BSVIEs.
Then the mapping $x\mapsto Y_{\Pi_N}^{t,x}(\cd)$ is differentiable with
the derivative terms $\nabla Y_{\Pi_N}^{t,x}(\cd)$ given by
\begin{align*}
\nabla Y^{t,x}_{\Pi_N}(s):=\left\{\begin{aligned}
&\nabla Y^{t,x}_{N}(s),\qq s\in[t_{N-1},T];\\
&\nabla Y^{t,x}_{N-1}(s),\qq s\in[t_{N-2},t_{N-1});\\
&\qq \cdot\cdot\cdot\cdot\cdot\cdot\\
&\nabla Y^{t,x}_{1}(s),\qq s\in[0,t_1),
\end{aligned}\right.\\
\nabla \cY^{t,x}_{\Pi_N}(s):=\left\{\begin{aligned}
&\nabla \cY^{t,x}_{N}(s),\qq s\in[t_{N-1},T];\\
&\nabla \cY^{t,x}_{N-1}(s),\qq s\in[t_{N-2},t_{N-1});\\
&\qq \cdot\cdot\cdot\cdot\cdot\cdot\\
&\nabla \cY^{t,x}_{1}(s),\qq s\in[0,t_1),
\end{aligned}\right.\end{align*}
where  for $k=1,2,...,N$,
$$\left\{\begin{aligned}
\nabla Y_{k}^{t,x}(s)&=h_x(X^{t,x}(T))+\int_s^T\[ g_x(r)+g_y(r)\nabla Y_{k}^{t,x}(r)+g_{y^\prime}(r) \nabla\cY_{k}^{t,x}(r) +g_z(r)\nabla Z_{k}^{t,x}(r)\]dr\\
&\q-\int_s^T\nabla Z_{k}^{t,x}(r)dW(r),\qq s\in[t_{k-1},T],\\
\nabla\cY_{k}^{t,x}(s)&=\int_{s\vee t_{k}}^T M(r-s)\nabla Y^{r,X^{t,x}(s)}_{\Pi_N}(r)dr,\qq s\in[t_{k-1},T].\end{aligned}\right.
$$
In the above, we used the equality $\nabla X^{t,x}(\cd)\equiv I$.
Then by the standard estimates of BSDEs, we have
\begin{align}
&|\nabla Y_{k}^{t,x}(t)|^2+\dbE_t\[\int_t^T |\nabla Z_{k}^{t,x}(s)|^2 ds\]
\les K+K\dbE_t\[\int_t^T |\nabla\cY_k^{t,x}(s)|^2 ds\]\nonumber\\
&\q\les   K+K\dbE_t\[\int_t^T \int_{s\vee t_{k}}^T |\nabla Y^{r,X^{t,x}(s)}_{\Pi_N}(r)|^2 drds\]\nn\\
&\q\les K+K\int_t^T \sup_{r,x\in[s,T]\times\dbR^n} |\nabla Y_{\Pi_N}^{r,x}(r)|^2 ds.\nn
%\label{nablaY}
\end{align}
%where $K$ is independent of $x$ and $N$.
It follows that
$$
|\nabla Y_{\Pi_N}^{t,x}(t)|^2
\les K+K\int_t^T \sup_{r,x\in[s,T]\times\dbR^n} |\nabla Y_{\Pi_N}^{r,x}(r)|^2 ds,\qq (t,x)\in[0,T]\times\dbR^n.
$$
Then by Gr\"{on}wall's inequality, we get
\bel{Step2-1}
|\nabla Y_{\Pi_N}^{t,x}(t)|^2\les K,\qq (t,x)\in[0,T]\times\dbR^n.
\ee

\ss

\emph{\textbf{Step 3.} Cauchy sequence with the norm  $\sup_{t,x}| Y_{\Pi_N}^{t,x}(t)|^2$.}

\ss

Bear in mind that  $\Pi_N: t_0=0,t_1={T\over N},......,t_{N-1}={(N-1)T\over N},t_N=T$, and we introduce another partition of $[0,T]$ by defining
$\Pi_{N^\prime}: t_0=0,t_1={T\over N^\prime},......,t_{N^\prime-1}={(N^\prime-1)T\over N^\prime},t_{N^\prime}=T$.
Suppose that $\Pi_N\subseteq \Pi_{N^\prime}$ and $ [t_{k^\prime},t_{k^\prime+1}) \subseteq [t_k,t_{k+1})$.
Then
$$\left\{\begin{aligned}
[Y_{k}^{t,x}-Y_{k^\prime}^{t,x}](s)&=\int_s^T \[g\big(r,X^{t,x}(r),Y_{k}^{t,x}(r),Z_{k}^{t,x}(r),\cY_{k}^{t,x}(r)\big)\\
&\qq-g\big(r,X^{t,x}(r),Y_{k^\prime}^{t,x}(r),Z_{k^\prime}^{t,x}(r),\cY_{k^\prime}^{t,x}(r)\big)\]dr\\
&\q-\int_s^T \big[Z_{k}^{t,x}(r)-Z_{k^\prime}^{t,x}(r)\big]dW(r), \q s\in  [t_{k^\prime},T],\\
 [\cY_{k}^{t,x}-\cY_{k^\prime}^{t,x}](s)&=\int_{s}^T M(r-s)[ Y^{r,X^{t,x}(s)}_{\Pi_N}-Y_{\Pi_{N^\prime}}^{t,X^{t,x}(s)}](r)dr,\q s\in [t_{k +1},T],\\
 [\cY_{k}^{t,x}-\cY_{k^\prime}^{t,x}](s)&=\int_{t_{k+1}}^T M(r-s)[ Y^{r,X^{t,x}(s)}_{\Pi_N}-Y_{\Pi_{N^\prime}}^{r,X^{t,x}(s)}](r)dr\\
 &\qq -\int^{t_{k+1}}_{s\vee t_{k^\prime+1}} M(r-s) Y^{r,X^{t,x}(s)}_{\Pi_{N^\prime}}(r)dr,\q s\in [t_{k^\prime},t_{k+1}]. \end{aligned}\right.
$$

By the standard estimates of BSDEs again, we have
\begin{align}
\dbE_t\[\sup_{s\in [t,T]}|Y_{k}^{t,x}(s)-Y_{k^\prime}^{t,x}(s)|^2\]
\les K \dbE_t\[\int_{t}^T|\cY_{k}^{t,x}(s)-\cY_{k^\prime}^{t,x}(s)|^2ds \],\nn
\end{align}
where $K$ is independent of $\Pi_N$ and $\Pi_{N^\prime}$.
From the definition of $[\cY_{k}^{t,x}-\cY_{k^\prime}^{t,x}](\cd)$,
we get
\begin{align}
\dbE_t\[\int_{t}^T|\cY_{k}^{t,x}(s)-\cY_{k^\prime}^{t,x}(s)|^2ds \] &\les {K\over N} \sup_{r,x}|Y^{r,x}_{\Pi_{N^\prime}}(r)|^2+K\int_{t}^T \sup_{x}| Y^{r,x}_{\Pi_N}(r)-Y_{\Pi_{N^\prime}}^{r,x}(r)|^2dr\nn\\
&\les {K\over N} +K\int_{t}^T \sup_{x}| Y^{r,x}_{\Pi_N}(r)-Y_{\Pi_{N^\prime}}^{r,x}(r)|^2dr,\nn
\end{align}
where the last inequality is due to the uniform estimate in Step 1. Thus,
$$
\dbE_t\[\sup_{s\in [t,T]}|Y_{k}^{t,x}(s)-Y_{k^\prime}^{t,x}(s)|^2\]
\les {K\over N} +K\int_{t}^T \sup_{x}| Y^{r,x}_{\Pi_N}(r)-Y_{\Pi_{N^\prime}}^{r,x}(r)|^2dr,
$$
which implies that
\begin{align}
&|Y_{k}^{t,x}(t)-Y_{k^\prime}^{t,x}(t)|^2\les {K\over N} +K\int_{t}^T \sup_{x}| Y^{r,x}_{\Pi_N}(r)-Y_{\Pi_{N^\prime}}^{r,x}(r)|^2dr,\nn
\end{align}
for $t\in  [t_{k^\prime},t_{k^\prime+1}) \subseteq [t_k,t_{k+1})$.
Then,
\bel{Stepp3-2}
\sup_{t,x}|Y_{\Pi_N}^{t,x}(t)-Y_{\Pi_{N^\prime}}^{t,x}(t)|^2\les {K\over N}.
\ee
This completes the proof.
\end{proof}

\subsection{Proof of Theorem \ref{thm:well-posdeness}}

We  are now ready to prove Theorem \ref{thm:well-posdeness}.

\ss

\emph{\textbf{Existence.}}
From \rf{Stepp3-2}, we have that $\{Y_{\Pi_N}^{\cd,\cd}(\cd)\}$ is a Cauchy sequence under the uniform norm.
Then the following deterministic  function is well-defined:
\bel{def-Ytxt}
 V(t,x):=\lim_{N\to\infty}Y_{\Pi_N}^{t,x}(t),\qq (t,x)\in[0,T]\times\dbR^n.
\ee
Moreover,
\bel{def-Ytxt1}
\lim_{N\to\infty}\sup_{ t,x}  |Y_{\Pi_N}^{t,x}(t)-V(t,x)|=0.
\ee
From \rf{Y-N-k},  we have
$$
Y_{\Pi_N}^{t^-_k,x}(t^-_k):=\lim_{t\to t_k^-}Y_{\Pi_N}^{t,x}(t)=\lim_{t\to t_k^-}Y_{k}^{t,x}(t)=Y_{k}^{t_k,x}(t_k),\qq
Y_{\Pi_N}^{t_k,x}(t_k)=Y_{k+1}^{t_k,x}(t_k).
$$
Thus, the mapping $t \mapsto Y_{\Pi_N}^{t,x}(t)$ has jumps at each $t_k$.
By the standard estimates of BSDEs, we have
\begin{align*}
|Y_{\Pi_N}^{t^-_k,x}(t^-_k)-Y_{\Pi_N}^{t_k,x}(t_k)|^2 =|Y_{k}^{t_k,x}(t_k)-Y_{k+1}^{t_k,x}(t_k)|^2
\les K\int_{t_k}^{T} |\cY_k^{t,x}(s) -\cY_{k+1}^{t,x}(s)|^2ds.
\end{align*}
Note that
\begin{align*}
&\cY_k^{t,x}(s) =\cY_{k+1}^{t,x}(s),\qq s\in[t_{k+1},T], \\
&\cY_k^{t,x}(s)-\cY_{k+1}^{t,x}(s)=\int_s^{t_{k+1}} M(r-s) Y^{r,X^{t,x}(s)}_{\Pi_N}(r)dr,\qq s\in[t_k, t_{k+1}].
\end{align*}
Thus,
$$
|Y_{\Pi_N}^{t^-_k,x}(t^-_k)-Y_{\Pi_N}^{t_k,x}(t_k)|^2 \les K\int_{t_k}^{t_{k+1}} |\cY_k^{t,x}(s) -\cY_{k+1}^{t,x}(s)|^2ds
\les {K\over N^2}.
$$
Now, for any $\e>0$, we can find an $N^*$ such that for $N>N^*$,
$$
 |Y_{\Pi_N}^{s,x}(s)-V(s,x)|\les {\e\over 4},\qq \forall s\in[0,T],\ x\in\mathbb{R}^n.
$$
Note that the mapping  $s\mapsto Y_{\Pi_N}^{s,x}(s)$ is continuous except the jump points.
Thus, for the fixed point $t$, we can find a $\d_N>0$ such that
$$
 |Y_{\Pi_N}^{s,x}(s)-Y_{\Pi_N}^{t,x}(t)|\les {\e\over 4}+{K\over N^2},\qq \forall s\in(t-\d_N,t+\d_N),\ x\in\mathbb{R}^n.
$$
Thus,
$$
|V(s,x)-V(t,x)|\les |Y_{\Pi_N}^{s,x}(s)-Y_{\Pi_N}^{t,x}(t)|+{\e\over 2}\les {3\e\over 4}+{K\over N^2},\qq \forall  s\in(t-\d_N,t+\d_N),\ x\in\mathbb{R}^n.
$$
Note that $K$ is independent of $s,x,t,N$. Taking $N>{4 K\over\e}$ and $\d=\d_N$, we get
\bel{V-tc}
|V(s,x)-V(t,x)|\les \e,\qq \forall  s\in(t-\d,t+\d).
\ee
Thus,  the mapping $t\mapsto  V(t,x)$ is continuous.
Moreover, by \rf{Step2-1} we have
$$
|Y_{\Pi_N}^{s,x_1}(s)-Y_{\Pi_N}^{s,x_2}(s)|\les \sup_{s,x} |\nabla Y_{\Pi_N}^{s,x}(s)||x_1-x_2|\les K|x-x^\prime|.
$$
Combining this with \rf{def-Ytxt}, we have
\bel{def-Vx}
|V(t,x)-V(t,y)|\les K|x-y|.
\ee
Since $s\mapsto X^{t,x}(s)$ is H\"{o}lder continuous, the above implies that
the mapping $s\mapsto \cY^{t,x}(s)$ is continuous, where
\bel{conv-step4-1}
\cY^{t,x}(s):=\int_{s}^T M(r-s) V(r,X^{t,x}(s))dr.
\ee
Now with the process $\cY^{t,x}(\cd)$ given by \rf{conv-step4-1},  the following BSDE admits a unique adapted solution
$(Y^{t,x}(\cd),Z^{t,x}(\cd))\in L^2_{\dbF}(\Omega; C([t,T];\dbR^n))\times L^2_{\dbF}([t,T];\dbR^n)$:
\bel{conv-step4-2}
\left\{\begin{aligned}
Y^{t,x}(s)&=h(X^{t,x}(T))+\int_s^T g\big(r,X^{t,x}(r),Y^{t,x}(r),Z^{t,x}(r),\cY^{t,x}(r)\big)dr-\int_s^TZ^{t,x}(r)dr,\\
\cY^{t,x}(s)&=\int_{s}^T M(r-s) V(r,X^{t,x}(s))dr.
\end{aligned}\right.
\ee
For the partition $\Pi_N\equiv\{t_k~|~0\les k\les N\}$, $(Y^{t,x}_1(\cd),Z^{t,x}_1(\cd))$ is given by
\bel{conv-step4-23}
\left\{\begin{aligned}
Y_{1}^{t,x}(s)&=h(X^{t,x}(T))+\int_s^T g\big(r,X^{t,x}(r),Y_1^{t,x}(r),Z_1^{t,x}(r),\cY_1^{t,x}(r)\big)dr-\int_s^TZ_1^{t,x}(r)dr,\\
\cY_1^{t,x}(s)&=\int_{s\vee t_1}^T M(r-s) Y^{r,X^{t,x}(s)}_{\Pi_N}(r)dr.
\end{aligned}\right.
\ee
Comparing BSDEs \rf{conv-step4-2} and \rf{conv-step4-23},  from \rf{def-Ytxt1}, by the stability estimates of BSDEs, we have
\bel{conv-step4-3}
\lim_{N\to\infty} \dbE\[\sup_{s\in[0,T]}|Y_{1}^{t,x}(s)-Y^{t,x}(s)|^2\]=0.
\ee
Next, for $s\in[t_{k-1},t_k)$, we consider the BSDE \rf{BSVIE-k}.
By the stability estimates of BSDEs again, we have
$$
|Y_{k}^{t,x}(s)-Y_1^{t,x}(s)|^2\les {K\over N}.
$$
Thus,
\bel{conv-step4-5}
|Y_{\Pi_N}^{t,x}(t)-Y_1^{t,x}(t)|^2\les {K\over N}.
\ee
Combining \rf{def-Ytxt1}, \rf{conv-step4-3}, and \rf{conv-step4-5}, we get
$$
V(t,x)=Y^{t,x}(t), \qq (t,x)\in[0,T]\times\dbR^n,
$$
which implies
$$
V(r,X^{t,x}(s))=Y^{r,\ti x}(y)|_{\ti x=X^{t,x}(s)}=Y^{r,X^{t,x}(s)}(r), \qq (t,x)\in[0,T]\times\dbR^n.
$$
Substituting this into \rf{conv-step4-2}, we get that $(Y^{t,x}(\cd),Z^{t,x}(\cd))$ satisfies the MBSDE \rf{MBSDE-rw}.
This completes the proof of existence of a solution.

\ss

\emph{\textbf{Uniqueness.}}
We next prove the uniqueness of solutions to MBSDE \rf{MBSDE-rw}.
Suppose that $(Y_1^{\cd,\cd}(\cd),Z_1^{\cd,\cd}(\cd))$ and  $(Y_2^{\cd,\cd}(\cd),Z_2^{\cd,\cd}(\cd))$
satisfy the MBSDE \rf{MBSDE-rw}. Then by the standard estimate of BSDEs, we have
\begin{align*}
\dbE_t\[\sup_{s\in[t,T]}|Y^{t,x}_1(s)-Y^{t,x}_2(s)|^2\]
&\les K\dbE_t\[\int_t^T \big|g\big(s,X^{t,x}(s),Y_2^{t,x}(s),Z_2^{t,x}(s),\cY_1^{t,x}(s)\big)\\
&\qq-g\big(s,X^{t,x}(s),Y_2^{t,x}(s),Z_2^{t,x}(s),\cY_2^{t,x}(s)\big)\big|^2 ds \]\\
&\les K\dbE_t\[\int_t^T \big|\cY_1^{t,x}(s)-\cY_2^{t,x}(s)|^2ds\]\nn\\
&\les K\int_t^T \sup_{(r,x)\in[s,T]\times\dbR^n}|Y_1^{r,x}(r)-Y_2^{r,x}(r)|^2ds.
\end{align*}
Thus,
$$
\sup_{(s,x)\in[t,T]\times\dbR^n}|Y_1^{s,x}(s)-Y_2^{s,x}(s)|^2\les K\int_t^T \sup_{(r,x)\in[s,T]\times\dbR^n}|Y_1^{r,x}(r)-Y_2^{r,x}(r)|^2ds,\qq t\in[0,T],
$$
which implies  that
$$
Y_1^{s,x}(s)=Y_2^{s,x}(s),\qq (s,x)\in[0,T]\times\dbR^n.
$$
Then by the definition of $\cY_i^{t,x}(\cd),i=1,2$, we have  $\cY_1^{t,x}(\cd)=\cY_2^{t,x}(\cd)$.
Further, by the uniqueness of solutions to BSDEs, we get $(Y_1^{\cd,\cd}(\cd),Z_1^{\cd,\cd}(\cd))=(Y_2^{\cd,\cd}(\cd),Z_2^{\cd,\cd}(\cd))$.
This completes the proof.

\subsection{Comparison theorem of MBSDEs}

Consider the following two one-dimensional MBSDEs:
\bel{2MBSDEs}\left\{\begin{aligned}
&dX^{t,x}(s)=\si(s) dW(s),\qq s\in[t,T],\\
&dY_i^{t,x}(s)=-g_i\big(s,X^{t,x}(s),Y_i^{t,x}(s),Z_i^{t,x}(s),\cM^{t,x}[Y_i](s)\big)ds\\
&\qq\qq\q+Z_i^{t,x}(s)dW(s),\qq s\in[t,T],\qq i=1,2,\\
&X^{t,x}(t)=x,\qq Y_i^{t,x}(T)=h_i(X^{t,x}(T)),\qq i=1,2.\end{aligned}\right.
\ee
Since the MBSDE is non-local, we introduce the following monotonicity condition, motivated by \cite{Wang-Yong2019,Wang-Yong-Zhang2022}, in order to establish its comparison theorem.

\begin{taggedassumption}{(H2)}\label{ass:H3}
The dimension of the backward equation is $1$, namely $m=1$ and $k=1$.
Suppose that  the mapping $y^\prime\mapsto g_1(t,x,y,z,y^\prime)$ is nondecreasing for any
$(t,x,y,z)\in [0,T]\times\dbR^n\times\dbR\times\dbR$, and
\begin{align*}
g_1(t,x,y,z,y^\prime)\ges g_2(t,x,y,z,y^\prime),\qq h_1(x)\ges h_2(x), \qq M(t)\ges 0,\\
\forall (t,x,y,z,y^\prime)\in [0,T]\times\dbR^n\times\dbR\times\dbR\times\dbR.
\end{align*}
\end{taggedassumption}

\begin{theorem}\label{thm:BSDE-Comparison}
Under {\rm\ref{ass:H1}} with $g(\cd)=g_i(\cd)$ and $h(\cd)=h_i(\cd)$, $i=1,2$, and \ref{ass:H3}, the solution $(Y^{\cd,\cd}_i(\cd),Z^{\cd,\cd}_i(\cd))$ to \rf{2MBSDEs}, $i=1,2$, satisfy
$$
Y_1^{t,x}(s)\ges Y_2^{t,x}(s),\qquad (t,s,x)\in\D[0,T]\times\dbR^n.
$$

\end{theorem}

\ms

To prove  Theorem \ref{thm:BSDE-Comparison}, we first establish the following lemma.
Consider the following BSVIEs:
\bel{BSDE-VIE}
\left\{\begin{aligned}
&Y_i^{t,x}(s)=h_i(X^{t,x}(T))+\int_s^T g_i\big(r,X^{t,x}(r),Y_i^{t,x}(r),Z_i^{t,x}(r),\cY_i^{t,x}(r)\big)dr\\
&\qq\qq-\int_s^T Z_i^{t,x}(r)dW(r),\qq s\in[t,T],\\
&\cY_i^{t,x}(s)=\int_s^T M_i(r-s)Y_i^{t,x}(r)dr-\int_s^T\cZ_i^{t,x}(s,r)dW(r),\qq s\in[t,T],
\end{aligned}\right.
\ee
with $i=1,2$.

\ms

\begin{lemma}\label{lem:comparison}
Under {\rm\ref{ass:H1}} with $g(\cd)=g_i(\cd)$ and $h(\cd)=h_i(\cd)$, $i=1,2$, and \ref{ass:H3},
for any given $(t,x)\in[0,T)\times\dbR^n$, the solution $(Y^{t,x}_i(\cd),Z^{t,x}_i(\cd))$ to \rf{BSDE-VIE}, $i=1,2$,  satisfy
\bel{lem:comparison-main1}
Y_1^{t,x}(s)\ges Y_2^{t,x}(s),\qq \forall (s,x)\in[t,T]\times\dbR^n.
\ee
\end{lemma}

\begin{proof}
Let $\ti Y^{t,x}_0(\cd)=Y_2^{t,x}(\cd)$, and for $k=1,2,...$,
$$
\left\{\begin{aligned}
&\ti Y_k^{t,x}(s)=h_1(X^{t,x}(T))+\int_s^T g_1\big(r,X^{t,x}(r),\ti Y_{k}^{t,x}(r),\ti Z_k^{t,x}(r),\ti\cY_{k-1}^{t,x}(r)\big)dr\\
&\qq\qq-\int_s^T \ti Z_k^{t,x}(r)dW(r),\qq s\in[t,T],\\
&\ti\cY_{k-1}^{t,x}(s)=\int_s^T M_1(r-s)\ti Y_{k-1}^{t,x}(r)dr-\int_s^T\ti\cZ_{k-1}^{t,x}(s,r)dW(r),\qq s\in[t,T].
\end{aligned}\right.
$$
Then by a standard iterative estimates of BSVIEs, we know that
\bel{lem:comparison-proof4}
\lim_{k\to\infty}\dbE\[\sup_{s\in[t,T]} |\ti Y_k^{t,x}(s)-Y_1^{t,x}(s)|^2\]=0.
\ee
Note that
$$
\ti\cY_0^{t,x}(s)=\cY_2^{t,x}(s),\qq s\in[t,T],
$$
and then
\begin{align}
&\ti Y_1^{t,x}(s)=h_1(X^{t,x}(T))+\int_s^T g_1\big(r,X^{t,x}(r),\ti Y_1^{t,x}(r),\ti Z_1^{t,x}(r),\cY_2^{t,x}(r)\big)dr\nn\\
&\qq\qq-\int_s^T \ti Z_1^{t,x}(r)dW(r),\label{lem:comparison-proof1}\\
&Y_2^{t,x}(s)=h_2(X^{t,x}(T))+\int_s^T g_2\big(r,X^{t,x}(r),Y_2^{t,x}(r),Z_2^{t,x}(r),\cY_2^{t,x}(r)\big)dr\nn\\
&\qq\qq-\int_s^T Z_2^{t,x}(r)dW(r),\label{lem:comparison-proof2}
\end{align}
Then with the fixed $\cY_2^{t,x}(\cd)$, the above can be regarded as two BSDEs with the unknown
$(\ti Y^{t,x}(\cd),\ti Z^{t,x}(\cd))$ and $(Y_2^{t,x}(\cd),Z_2^{t,x}(\cd))$.
Then by the standard comparison theorem of BSDEs, we know that
$$
\ti Y_1^{t,x}(s)\ges Y_2^{t,x}(s),\qq s\in[t,T].
$$
Then by $M(\cd)\ges 0$, we have $\ti\cY^{t,x}_1(\cd)\ges \cY_2^{t,x}(\cd)$,
which, together with the monotonicity of the mapping $y^\prime\mapsto g_1(t,x,y,z,y^\prime)$, implies
\begin{align*}
 g_1(r,x,y,z,\ti\cY_1^{t,x}(r))\ges  g_1(r,x,y,z,\cY_2^{t,x}(r))\ges g_2(r,x,y,z,\cY_2^{t,x}(r)),\\
 \forall (r,x,y,z)\in [t,T]\times\dbR\times\dbR\times\dbR.
\end{align*}
Note that
\begin{align}
&\ti Y_2^{t,x}(s)=h_1(X^{t,x}(T))+\int_s^T g_1\big(r,X^{t,x}(r),\ti Y_2^{t,x}(r),\ti Z_2^{t,x}(r),\ti\cY_1^{t,x}(r)\big)dr\nn\\
&\qq\qq-\int_s^T \ti Z_2^{t,x}(r)dW(r).\label{lem:comparison-proof3}
\end{align}
Applying the comparison theorem of BSDEs for \rf{lem:comparison-proof1}--\rf{lem:comparison-proof3} again,
we get
$$
\ti Y_2^{t,x}(s)\ges Y_2^{t,x}(s),\qq s\in[t,T].
$$
Repeating the above arguments, we get
$$
\ti Y_k^{t,x}(s)\ges Y_2^{t,x}(s),\qq s\in[t,T],
$$
for any $k\ges 1$. Combining the above with \rf{lem:comparison-proof4}, we get
\rf{lem:comparison-main1} immediately.
\end{proof}

\begin{remark}
Lemma \ref{lem:comparison} can be proved by directly following Wang--Yong \cite{Wang-Yong2019}.
For completeness, a sketch of the proof is provided above.
\end{remark}

\noindent
\emph{\textbf{Proof of Theorem \ref{thm:BSDE-Comparison}.}}
We consider the approximating equation of as that in Subsection \ref{subsec:approximate}.
Let $\Pi_N\equiv\{t_k~|~0\les k\les N\}$ be a partition of $[0,T]$ with
$t_k={k\over N}T$.
On $[t_{N-1},t_N]$, let
$$
\left\{\begin{aligned}
&Y_{i,N}^{t,x}(s)=h_i(X^{t,x}(T))+\int_s^T g_i\big(r,X^{t,x}(r),Y_{i,N}^{t,x}(r),Z_{i,N}^{t,x}(r),\cY_{i,N}^{t,x}(r)\big)dr-\int_s^T Z_{i,N}^{t,x}(r)dW(r),\\
&\cY_{i,N}^{t,x}(s)=\int_s^T M(r-s)Y_{i,N}^{t,x}(r)dr-\int_s^T\cZ_{i,N}^{t,x}(s,r)dW(r).
\end{aligned}\right.
$$
By Lemma \ref{lem:comparison}, we have
\bel{thm:BSDE-Comparison-p1}
Y_{1,N}^{t,x}(s)\ges Y_{2,N}^{t,x}(s),\qq s\in[t_{N-1}\vee t,T],
\ee
for any $(t,x)\in[0,T)\times\dbR^n$. Then on $[t_{N-2},T]$, we have
\bel{thm:BSDE-Comparison-p2}
\left\{\begin{aligned}
Y_{i,N-1}^{t,x}(s)&=h_i(X^{t,x}(T))+\int_s^T g_i\big(r,X^{t,x}(r),Y_{i,N-1}^{t,x}(r),Z_{i,N-1}^{t,x}(r),\cY_{i,N-1}^{t,x}(r)\big)dr\\
&\q-\int_s^TZ_{i,N-1}^{t,x}(r)dr,\qq s\in[t_{N-2},T],\\
\cY_{i,N-1}^{t,x}(s)&=\int_{s\vee t_{N-1}}^T M(r-s)Y^{r,X^{t,x}(s)}_{i,N}(r)dr,\qq s\in[t_{N-2},T].\end{aligned}\right.
\ee
Note that $M(\cd)\ges 0$, from \rf{thm:BSDE-Comparison-p1} we get
$$
\cY_{1,N-1}^{t,x}(s)\ges \cY_{2,N-1}^{t,x}(s),\qq s\in[t_{N-2},T],
$$
which, together with \ref{ass:H3}, yields that
\begin{align}
g_1\big(r,x,y,z,\cY_{i,N-1}^{t,x}(r)\big)\ges g_2\big(r,x,y,z,\cY_{2,N-1}^{t,x}(r)\big),\nn\\
 \forall (r,x,y,z)\in [t_{N-2},T]\times\dbR\times\dbR\times\dbR.\label{thm:BSDE-Comparison-p3}
\end{align}
With the fixed $\cY_{i,N-1}^{t,x}(\cd)$, \rf{thm:BSDE-Comparison-p2} can be regarded as a standard BSDE.
Then by the standard comparison theorem of BSDEs, from \rf{thm:BSDE-Comparison-p3} we get
\bel{thm:BSDE-Comparison-p4}
Y_{1,N-1}^{t,x}(s)\ges Y_{2,N-1}^{t,x}(s),\qq s\in[t_{N-2},T].
\ee
Define
$$
Y^{s,x}_{i,\Pi_N}(s):=\left\{\begin{aligned}
&Y^{s,x}_{i,N}(s),\qq s\in[t_{N-1},T],\q x\in\dbR^n,\\
&Y^{s,x}_{i,N-1}(s),\qq s\in[t_{N-2},t_{N-1}),\q x\in\dbR^n.
\end{aligned}\right.
$$
Then by \rf{thm:BSDE-Comparison-p1} and \rf{thm:BSDE-Comparison-p4}, we have
\bel{thm:BSDE-Comparison-p8}
Y^{s,x}_{1,\Pi_N}(s)\ges Y^{s,x}_{2,\Pi_N}(s), \qq (s,x)\in[t_{N-2},T]\times\dbR^n.
\ee
On the general time horizon $[t_{k-1},T]$, we consider
$$
\left\{\begin{aligned}
Y_{k}^{t,x}(s)&=h(X^{t,x}(T))+\int_s^T g\big(r,X^{t,x}(r),Y_k^{t,x}(r),Z_k^{t,x}(r),\cY_k^{t,x}(r)\big)dr-\int_s^TZ_{k}^{t,x}(r)dr,\\
\cY_{k}^{t,x}(s)&=\int_{s\vee t_{k}}^T M(r-s) Y^{r,X^{t,x}(s)}_{\Pi_N}(r)dr.
\end{aligned}\right.
$$
With \rf{thm:BSDE-Comparison-p8}, we have
%\bel{thm:BSDE-Comparison-p9}
$$
Y_{1,k}^{t,x}(s)\ges Y_{2,k}^{t,x}(s),\qq s\in[t_{k-1},T].
$$
By induction, we have
$$
Y_{1,1}^{t,x}(s)\ges Y_{2,1}^{t,x}(s),\qq s\in[0,T].
$$
Combining the above with \rf{conv-step4-3}, we have
$$
Y_{1}^{t,x}(s)\ges Y_{2}^{t,x}(s),\qq s\in[0,T],
$$
which completes the proof.
$\hfill\blacksquare$

\section{Classical solutions}\label{sec:6}

In this section, we shall establish the Feynman--Kac formula in the sense of classical solutions.
Under \ref{ass:H1}, from Theorem \ref{thm:well-posdeness}, the MBSDE   \rf{MBSDE-rw} admits a unique adapted solution $(Y^{\cd,\cd}(\cd),Z^{\cd,\cd}(\cd)) $.
Define
\bel{u-Y}
u(t,x)=Y^{t,x}(t),\qq (t,x)\in[0,T]\times\dbR^n.
\ee

\begin{proposition}\label{prop:VT}
Suppose that the function $u(\cd,\cd)$ defined by \rf{u-Y} belongs to $ C^{0,2}([0,T]\times\dbR^n;\dbR^m)$.
Then it is a classical solution to the MPE \rf{H-PDE-rw}--\rf{def-M-operator-rw}.
\end{proposition}

\begin{proof}
By the flow property of SDEs, we have that $X^{t,x}(r)=X^{s,X^{t,x}(s)}(r)$ for any $(t,x)\in[0,T]\times\dbR^n$ and $t\les s\les r\les T$.
Then by the uniqueness of solutions to MBSDEs, we have
$$
Y^{t,x}(r)=Y^{s,X^{t,x}(s)}(r),\qq (t,x)\in[0,T]\times\dbR^n,\q t\les s\les r\les T.
$$
In particular, taking $r=s$ in the above, we have
$$
Y^{t,x}(s)=Y^{s,X^{t,x}(s)}(s)=u(s,X^{t,x}(s)),\qq (t,x)\in[0,T]\times\dbR^n,\q t\les s\les T.
$$
We next show that
\bel{prop:VT-proof11}
Z^{t,x}(s)=\partial_{x}u(s,X^{t,x}(s))\si(s),\qq (t,x)\in[0,T]\times\dbR^n,\q t\les s\les T.
\ee
Fix $\d>0$ and let $s=s_0<\cdots<s_n=T$ with $\D s_i := s_i-s_{i-1}\les \d$.
Denote
$$
Y_n^{t,x}(s):=u(s_{i+1},X^{t,x}(s)),\qq
Z_n^{t,x}(s):=\partial_{x}u(s_{i+1},X^{t,x}(s))\si(s),\qq s_{i}\les s\les s_{i+1}.
$$
Note that $u(\cd,\cd)$ is a deterministic function and then $u(s_{i+1},X^{t,x}(s))$ is $\cF_s$-measurable.
By applying It\^{o}'s formula, we have
\begin{align*}
dY_n^{t,x}(s)&=du(s_{i+1},X^{t,x}(s))\\
&={1\over 2}\tr\big[\partial_{xx}^2 u(s_{i+1},X^{t,x}(s))\si(s)\si(s)^\top\big]ds+\partial_{x}u(s_{i+1},X^{t,x}(s))\si(s)dW(s)\\
&={1\over 2}\tr\big[\partial_{xx}^2 u(s_{i+1},X^{t,x}(s))\si(s)\si(s)^\top\big]ds+Z^{t,x}_n(s)dW(s),\qq s\in[s_i,s_{i+1}].
\end{align*}
Recall that
$$
dY^{t,x}(s)=-g\big(s,X^{t,x}(s),Y^{t,x}(s),Z^{t,x}(s),\cM^{t,x}[Y](s)\big)ds+Z^{t,x}(s)dW(s),\qq s\in[t,T].
$$
Denote $\D Y^{t,x}_n(s)=Y_n^{t,x}(s)-Y^{t,x}(s)$ and $\D Z^{t,x}_n(s)=Z_n^{t,x}(s)-Z^{t,x}(s)$. Then
\bel{prop:VT-proof1}
\left\{\begin{aligned}
&d\D Y^{t,x}_n(s)=\h g\big(s,X^{t,x}(s),Y^{t,x}(s),Z_n^{t,x}(s)-\D Z_n^{t,x}(s),\cM^{t,x}[Y](s)\big)ds\\
&\qq\qq\qq +\D Z^{t,x}_n(s) dW(s),\qq s\in[s_i,s_{i+1}],\\
&\D Y^{t,x}_n(s_{i+1})=0,
\end{aligned}\right.
\ee
where
$$
\h g(s,x,y,z,y^\prime)={1\over 2}\tr\big[\partial_{xx}^2 u(s_{i+1},x)\si(s)\si(s)^\top\big]+g(s,x,y,z,y^\prime).
$$
By the standard estimates of BSDEs, we get
\begin{align}
&\dbE\[\sup_{s\in[s_i,s_{i+1}]}|\D Y^{t,x}_n(s)|^2+\int_{s_i}^{s_{i+1}}|\D Z^{t,x}_n(s)|^2ds\]\nn\\
&\q\les K\dbE\[ \(\int_{s_i}^{s_{i+1}}\big|\h g\big(s,X^{t,x}(s),Y^{t,x}(s),Z_n^{t,x}(s),\cM^{t,x}[Y](s)\big)\big | ds\)^2\]\nn\\
&\q\les K \d\D s_{i+1}. \nn
\end{align}
Thus,
$$
\dbE\[\sum_{i=0}^{n-1} \int_{s_i}^{s_{i+1}}\big| Z^{t,x}(s)-\partial_{x}u(s_{i+1},X^{t,x}(s))\si(s)\big|^2ds\]\les K\d.
$$
Letting $\d\to 0$, by the continuity of  $\partial_{x}u(\cd,\cd)$, we get \rf{prop:VT-proof11}.

\ss
Taking $s_i=t$ and $s_{i+ 1}= t+\d$ in \rf{prop:VT-proof1},  we have
\begin{align}
&u(t+\d,x)-u(t,x)=\D Y^{t,x}_n(t)\nn\\
&\q=-\dbE\[\int_{t}^{t+\d}\h g\big(s,X^{t,x}(s),Y^{t,x}(s),Z_n^{t,x}(s)-\D Z_n^{t,x}(s),\cM^{t,x}[Y](s)\big)ds\]\nn\\
&\q=-\dbE\[\int_{t}^{t+\d}\bar g\big(s,X^{t,x}(s),Y^{t,x}(s),Z_n^{t,x}(s),\cM^{t,x}[Y](s)\big)ds\]-R(\d),\label{prop:VT-proof2}
\end{align}
where
\begin{align}
&\bar g(s,x,y,z,y^\prime)={1\over 2}\tr\big[\partial_{xx}^2 u(s,x)\si(s)\si(s)^\top\big]+g(s,x,y,z,y^\prime),\nn\\
& R(\d):=\dbE\[\int_t^{t+\d}\(\h g\big(s,X^{t,x}(s),Y^{t,x}(s),Z_n^{t,x}(s)-\D Z_n^{t,x}(s),\cM^{t,x}[Y](s)\big)\nn\\
&\qq\qq -\h g\big(s,X^{t,x}(s),Y^{t,x}(s),Z_n^{t,x}(s),\cM^{t,x}[Y](s)\big)\)ds\nn\\
&\qq\qq+\int_t^{t+\d}\(\h g\big(s,X^{t,x}(s),Y^{t,x}(s),Z_n^{t,x}(s),\cM^{t,x}[Y](s)\big)\nn\\
&\qq\qq-\bar g\big(s,X^{t,x}(s),Y^{t,x}(s),Z_n^{t,x}(s),\cM^{t,x}[Y](s)\)ds\].\nn
\end{align}
Note that
$$
|R(\d)|^2\les K \d\dbE\[\int_t^{t+\d} |\D Z_n^{t,x}(s)|^2ds\]+o(\d^2)\les K\d^3+o(\d^2)=o(\d^2).
$$
Then from \rf{prop:VT-proof2}, we have
\begin{align}
{u(t+\d,x)-u(t,x)\over \d}&=-{1\over \d}\dbE\[\int_{t}^{t+\d}\bar g\big(s,X^{t,x}(s),Y^{t,x}(s),Z_n^{t,x}(s),\cM^{t,x}[Y](s)\big)ds\]-o(\d).\label{prop:VT-proof3}
\end{align}
Recall that
$$
\cM^{t,x}[Y](s)=\int_s^T M(r-s) Y^{r,X^{t,x}(s)}(r)dr=\int_s^T  M(r-s) u(r,X^{t,x}(s))dr.
$$
Then  by the desired regularity of $u(\cd,\cd)$, taking $\d\to 0$ in \rf{prop:VT-proof3},
we see that $\partial_t u(\cd,\cd)$ exists and
\begin{align*}
\partial_t u(t,x)&=-\bar g\big(t,x,u(t,x),\partial_x u(t,x)\si(t),\cM[u](t,x)\big)\\
%\
&=-{1\over 2}\tr[\partial_{xx}^2 u(t,x)\si(t)\si(t)^\top]-g\big(t,x,u(t,x),\partial_x u(t,x)\si(t),\cM[u](t,x)\big).
\end{align*}
This completes the proof.
\end{proof}

To establish the $C^{0,2}$-regularity of $u(\cd,\cd)$, we introduce the following assumption.

\begin{taggedassumption}{(H3)}\label{ass:H4}
The mappings $x\mapsto h(x)$ and  $(x,y,z,y^\prime)\mapsto g(t,x,y,z,y^\prime)$ are smooth enough and the related derivatives  are all bounded.
\end{taggedassumption}

\begin{theorem}\label{Thm-FK-classical}
Under \ref{ass:H1} and \ref{ass:H4},
the function $u(\cd,\cd)$ defined by \rf{u-Y}  is a classical solution to the MPE \rf{H-PDE-rw}--\rf{def-M-operator-rw}.
\end{theorem}

\begin{proof}
By Proposition \ref{prop:VT}, it suffices to verify that $u(\cd,\cd)$ belongs to $ C^{0,2}([0,T]\times\dbR^n;\dbR^m)$.

\ss
\textbf{\emph{Step 1.}} The function $u(\cd,\cd)$ is differentiable with respect to $x$, and
\bel{Thm-FK-classical-p1}
\partial_x u(t,x)=\nabla Y^{t,x}(t),
\ee
where
$$\left\{\begin{aligned}
\nabla Y^{t,x}(s)&=h_x(X^{t,x}(T))+\int_s^T\[ g_x(r)+g_y(r)\nabla Y^{t,x}(r)+g_z(r)\nabla Z^{t,x}(r)+g_{y^\prime}(r) \nabla\cY^{t,x}(r)\]dr\\
&\q -\int_s^T\nabla Z^{t,x}(r)dW(r),\qq s\in[t,T],\\
\nabla\cY^{t,x}(s)&=\int_{s}^T M(r-s)\nabla Y^{r,X^{t,x}(s)}(r)dr,\qq s\in[t,T].\end{aligned}\right.
$$
We now prove the above claim.
For any fixed $\eta\in\dbR^n$,  denote
\begin{align*}
&\D  Y^{t,x}(s):=Y^{t,x+\eta}(s)-Y^{t,x}(s)-\nabla Y^{t,x}(s)\eta,\\
& \D  Z^{t,x}(s):=Z^{t,x+\eta}(s)-Z^{t,x}(s)-\nabla Z^{t,x}(s)\eta,\\
& \D \cY  ^{t,x}(s):=\cY^{t,x+\eta}(s)-\cY^{t,x}(s)-\nabla \cY^{t,x}(s)\eta.
\end{align*}
Then
$$\left\{\begin{aligned}
\D Y^{t,x}(s)&= h( X^{t,x+\eta}(T))-h( X^{t,x}(T))-h_x(X^{t,x}(T))\eta\\
 &\q+\int_s^T \ti g(r, \eta, \D  Y^{t,x}(r),\D  Z^{t,x}(r),\D \cY^{t,x}(r))dr-\int_s^T\D Z^{t,x}(r)dW(r),\qq s\in[t,T],\\
\D\cY^{t,x}(s)&=\int_{s}^T M(r-s)\D Y^{r,X^{t,x}(s)}(r)dr,\qq s\in[t,T],\end{aligned}\right.
$$
where
\begin{align*}
\ti g(r,\eta, y,z,y^\prime)&=g\big(r,X^{t,x+\eta}(r),Y^{t,x}(r)+\nabla Y^{t,x}(r)\eta+y,\, Z^{t,x}(r)+\nabla Z^{t,x}(r)\eta+z,\\
&\qq\cY^{t,x}(r)+\nabla \cY^{t,x}(r)\eta+y^\prime \big)-g\big(r,X^{t,x}(r),Y^{t,x}(r), Z^{t,x}(r),\cY^{t,x}(r) \big)\\
&\q -\big[ g_x(r)\eta +g_y(r)\nabla Y^{t,x}(r)\eta+g_{y^\prime}(r) \nabla\cY^{t,x}(r)\eta+g_z(r)\nabla Z^{t,x}(r)\eta\big].
\end{align*}
By the estimate in Theorem \ref{thm:well-posdeness}, we have
\begin{align}
&|\D Y^{t,x}(t)|^p\les K_p \dbE_t\[|h( X^{t,x+\eta}(T))-h( X^{t,x}(T))-h_x(X^{t,x}(T))\eta|^p+\(\int_t^T |\ti g(r,\eta, 0,0,0)|dr\)^p\]\nn\\
&\q\les K_p\dbE_t\[ \int_0^1\big| h_x\big( X^{t,x}(T)+\a\eta \big)- h_x\big( X^{t,x}(T) \big)\big|^p d\a |\eta|^p\]\nn\\
&\qq\q +K_p \dbE_t\[\int_t^T\( \Big| \int_0^1 g_x\big( r, X^{t,x}(r)+\a\eta,Y^{t,x}(r), Z^{t,x}(r),\cY^{t,x}(r) \big)d\a-g_x(r)\Big| |\eta|\nn\\
&\qq\q+ \Big| \int_0^1 g_x\big( r, X^{t,x}(r),Y^{t,x}(r)+\a\nabla Y^{t,x}(r)\eta , Z^{t,x}(r),\cY^{t,x}(r) \big)d\a-g_y(r)\Big|
|\nabla Y^{t,x}(r)\eta|\nn\\
&\qq\q +\Big| \int_0^1 g_z\big( r, X^{t,x}(r),Y^{t,x}(r), Z^{t,x}(r)+\a\nabla Z^{t,x}(r)\eta,\cY^{t,x}(r) \big)d\a-g_z(r)\Big|
|\nabla Z^{t,x}(r)\eta|\nn\\
&\qq\q +\Big| \int_0^1 g_{y^\prime}\big( r, X^{t,x}(r),Y^{t,x}(r), Z^{t,x}(r),\cY^{t,x}(r) +\a\nabla \cY^{t,x}(r)\eta\big)d\a-g_{y^\prime}(r)\Big|
|\nabla \cY^{t,x}(r)\eta|\)dr\]^p\nn\\
&\q\les K_p\dbE_t\[1+\int_t^T\( |\nabla Y^{t,x}(r)|^{2p}+ |\nabla Z^{t,x}(r)|^{2p}+ |\nabla \cY^{t,x}(r)|^{2p} \)dr \]|\eta|^{2p}\nn
\end{align}
By Theorem \ref{thm:well-posdeness}, we have
$$
\dbE_t\int_t^T\( |\nabla Y^{t,x}(r)|^{2p}+ |\nabla Z^{t,x}(r)|^{2p}+ |\nabla \cY^{t,x}(r)|^{2p} \)dr \les K_{2p}.
$$
Combining the above two inequalities, we have
$$
|\D Y^{t,x}(t)|^p\les K_p |\eta|^{2p}.
$$
The claim \rf{Thm-FK-classical-p1} is proved.

\ss
\textbf{\emph{Step 2.}}
By similar arguments as in Step 1, we can obtain that $u(\cd,\cd)$ is twice differentiable with respect to $x$, and
$$
\partial_{xx} u(t,x)=\nabla^2 Y^{t,x}(t),
$$
where
$$\left\{\begin{aligned}
\nabla^2 Y^{t,x}(s)&=h_{xx}(X^{t,x}(T))+\int_s^T\[ g_{xx}(r)+g_{xy}(r)\nabla Y^{t,x}(r)+g_{xz}(r)\nabla Z^{t,x}(r)\\
&\q+g_{xy^\prime}(r) \nabla\cY^{t,x}(r) +\( g_{yx}(r)+g_{yy}(r)\nabla Y^{t,x}(r)+g_{yz}(r)\nabla Z^{t,x}(r) \\
&\q+g_{yy^\prime}(r) \nabla\cY^{t,x}(r) \)\nabla Y^{t,x}(r)+ g_{y}(r) \nabla^2 Y^{t,x}(r)+ \( g_{zx}(r)+g_{zy}(r)\nabla Y^{t,x}(r)\\
&\q+g_{zz}(r)\nabla Z^{t,x}(r)+g_{zy^\prime}(r) \nabla\cY^{t,x}(r) \)\nabla Z^{t,x}(r)+ g_{z}(r) \nabla^2 Z^{t,x}(r) \\
&\q+ \( g_{y^\prime x}(r)+g_{y^\prime y}(r)\nabla Y^{t,x}(r)+g_{y^\prime z}(r)\nabla Z^{t,x}(r)+g_{y^\prime y^\prime}(r) \nabla\cY^{t,x}(r) \)\nabla \cY^{t,x}(r)\\
&\q+g_{y^\prime}(r) \nabla^2\cY^{t,x}(r)\]dr-\int_s^T\nabla^2 Z^{t,x}(r)dW(r),\qq s\in[t,T],\\
\nabla^2\cY^{t,x}(s)&=\int_{s}^T M(r-s)\nabla^2 Y^{r,X^{t,x}(s)}(r)dr,\qq s\in[t,T].\end{aligned}\right.
$$

\ss

\textbf{\emph{Step 3.}} We next show the continuity of $u(\cd,\cd)$, $\partial_x u(\cd,\cd)$, and $\partial_{xx} u(\cd,\cd)$.
For any $\eta\in\dbR^n$, by the same arguments as in Step 1, we have
$$
|u(t,x+\eta)-u(t,x)|=|Y^{t,x+\eta}(t)-Y^{t,x}(t)|\les K|\eta|,
$$
which means that the mapping $x\mapsto u(t,x)$ is Lipschitz continuous.
Similarly, we also have that $\partial_x u(\cd,\cd)$ and $\partial_{xx} u(\cd,\cd)$ are Lipschitz continuous
with respect to $x$. Next, for any $\d>0$, we have
\begin{align}
&u(t+\d,x)-u(t,x)=u(t+\d,x)-u(t+\d,X^{t,x}(t+\d))+u(t+\d,X^{t,x}(t+\d))-u(t,x)\nn\\
&\q=u(t+\d,x)-u(t+\d,X^{t,x}(t+\d))+Y^{t+\d,X^{t,x}(t+\d)}(t+\d)-Y^{t,x}(t)\nn\\
&\q=u(t+\d,x)-u(t+\d,X^{t,x}(t+\d))+Y^{t,x}(t+\d)-Y^{t,x}(t)\nn\\
&\q=u(t+\d,x)-u(t+\d,X^{t,x}(t+\d))+\int_t^{t+\d}Z^{t,x}(s)dW(s)\nn\\
&\qq-\int_t^{t+\d} g\big(s,X^{t,x}(s),Y^{t,x}(s),Z^{t,x}(s),\cM^{t,x}[Y](s)\big) ds.\nn
\end{align}
Thus,
\begin{align*}
&|u(t+\d,x)-u(t,x)|\\
&\q\les \dbE\[ |u(t+\d,x)-u(t+\d,X^{t,x}(t+\d))|+\int_t^{t+\d} \big|g\big(s,X^{t,x}(s),Y^{t,x}(s),Z^{t,x}(s),\cY^{t,x}(s)\big)\big| ds\]\\
&\q\les K\dbE\[|X^{t,x}(t+\d)-x|+\int_t^{t+\d}\big(1+|X^{t,x}(s)|+|Y^{t,x}(s)|+|Z^{t,x}(s)|+|\cY^{t,x}(s)|\big) ds\]\\
&\q \les K\d^{1\over 2}+\dbE\[\int_t^{t+\d}\big(1+|X^{t,x}(s)|^2+|Y^{t,x}(s)|^2+|Z^{t,x}(s)|^2+|\cY^{t,x}(s)|^2\big) ds\]^{1\over 2}\d^{1\over 2}\\
&\q\les K\d^{1\over2}.
\end{align*}
By the same arguments as the above, we have
$$
|\partial_x u(t+\d,x)-\partial_x u(t,x)|+|\partial_{xx}u(t+\d,x)-\partial_{xx}u(t,x)|\les  K\d^{1\over2}.
$$
This completes the proof.
\end{proof}

\section{Viscosity solutions}\label{sec:7}
In Theorem \ref{Thm-FK-classical}, we proved that under \ref{ass:H1} and \ref{ass:H4}, the function $u(\cd,\cd)$ defined by \rf{u-Y} belongs to $C^{1,2}([0,T]\times\dbR^n;\dbR^m)$ and is therefore a classical solution to the MPE \rf{H-PDE-rw}--\rf{def-M-operator-rw}. However, such regularity of $u(\cd,\cd)$ does not hold in general when the smoothness assumption \ref{ass:H4} on the coefficients is not satisfied. Indeed, it follows from \rf{V-tc}--\rf{def-Vx} that, under the standard Lipschitz continuity assumption \ref{ass:H1}, the function $u(\cd,\cd)$ defined by \rf{u-Y} is merely continuous. Thus, in this section, we establish the connection between the MBSDE \rf{MBSDE-rw}--\rf{def-F-rw} and the MPE \rf{H-PDE-rw}--\rf{def-M-operator-rw} under \ref{ass:H1}, using viscosity solution theory.

\ss
For notational simplicity, we define the operator $\cL[\,\cd\,]$ associated with the MPE \rf{H-PDE-rw}--\rf{def-M-operator-rw} as follows:
\begin{align*}
\cL[u](t,x)&:=\partial_tu(t,x)+{1\over 2}\tr\big[\partial_{xx}^2 u(t,x)\si(t)\si(t)^\top\big]\nn\\
&\q+g\big(t,x,u(t,x),\partial_x u(t,x)\si(t),\cM[u](t,x)\big),\qq \forall u(\cd,\cd)\in C^{1,2},
\end{align*}
where the memory operator $\cM[\,\cd\,]$ is defined by \rf{def-M-operator-rw}.

\begin{definition}
A continuous function $u(\cd,\cd)$ is called a {\it viscosity subsolution} (resp.\ {\it supersolution}) to the MPE \rf{H-PDE-rw}--\rf{def-M-operator-rw} if
$$ u(T,x)\les h(x), \q \forall x\in\dbR^n,\qq
(\hbox{resp.~~} u(T,x)\ges h(x), \q \forall x\in\dbR^n),
$$
and if, for any smooth test function $\f(\cd,\cd)$ satisfying $\f(\cd,\cd) \ges u(\cd,\cd)$
(resp.\ $\f(\cd,\cd) \les u(\cd,\cd)$),
whenever $\f(t, x) = u(t,x)$ at some $( t,   x)\in[0,T]\times\dbR^n$,
we have
$$
\cL[\f]( t, x) \ges 0,
\hbox{~~(resp.~~}  \cL [\f]( t, x) \les 0 \hbox{)}.
$$
Further, if $u(\cd,\cd)$ is both a viscosity subsolution and a viscosity supersolution,  then it is called a {\it viscosity solution} to the MPE \rf{H-PDE-rw}--\rf{def-M-operator-rw}.
 \end{definition}

The key in the viscosity solution theory is the so-called {\it comparison theorem}. Since the  MPE \rf{H-PDE-rw}--\rf{def-M-operator-rw} is non-local, similar to Theorem \ref{thm:BSDE-Comparison}, we introduce the following monotonicity condition.

\begin{taggedassumption}{(H4)}\label{ass:H5}
The MPE \rf{H-PDE-rw}--\rf{def-M-operator-rw} is $1$-dimensional, namely $m=1$ and $k=1$.
Suppose that $M(\cd)\ges 0$ and the mapping $y^\prime\mapsto g_1(t,x,y,z,y^\prime)$ is nondecreasing for any
$(t,x,y,z)\in [0,T]\times\dbR^n\times\dbR\times\dbR$.
\end{taggedassumption}

We first have the consistency result.

\begin{lemma}\label{lem-consistency}
Suppose that {\rm\ref{ass:H5}} holds and  $u(\cd,\cd)$ is smooth.
Then $u(\cd,\cd)$ is a viscosity subsolution to the MPE \rf{H-PDE-rw}--\rf{def-M-operator-rw} if and only if
$$
\cL[u](t,x)\ges 0, \qq u(T,x)\les h(x).
$$
\end{lemma}

\begin{proof}
The ``only if'' part is obvious by simply taking the test function $\f(\cd,\cd) = u(\cd,\cd)$.
We now prove the ``if'' part. Fix a test function $\f(\cd,\cd)$ such that
$\f(\cd,\cd)\ges u(\cd,\cd)$ and $\f(t,x)=u(t,x)$ at the fixed point $(t,x)$.
Define $\tilde{\varphi}=\varphi-u$. Then
\begin{align*}
\partial_t\ti\f(t,x)=0,\qq \partial_x\ti\f(t,x)=0,\qq \partial_{xx}\ti\f(t,x)\ges 0.
\end{align*}
Thus, at the fixed point $(t,x)$, we have
\begin{align*}
\cL[\ti\f](t,x)&={1\over 2}\tr\big[\partial_{xx}^2 \ti\f(t,x)\si(t)\si(t)^\top\big]+g\big(t,x,\f(t,x),\partial_x \f(t,x)\si(t),\cM[\f](t,x)\big)\\
&\q-g\big(t,x,u(t,x),\partial_x u(t,x)\si(t),\cM[u](t,x)\big)\\
&={1\over 2}\tr\big[\partial_{xx}^2 \ti\f(t,x)\si(t)\si(t)^\top\big]+g\big(t,x,\f(t,x),\partial_x \f(t,x)\si(t),\cM[\f](t,x)\big)\\
&\q-g\big(t,x,\f(t,x),\partial_x \f(t,x)\si(t),\cM[u](t,x)\big).
\end{align*}
From $M(\cd)\ges 0$ and $\f(\cd,\cd)\ges u(\cd,\cd)$, it follows that
$$
\cM[\f](\cd,\cd)\ges \cM[u](\cd,\cd).
$$
Since $g(\cd)$ is nondecreasing in $y^\prime$, we have
\begin{align*}
g\big(t,x,\f(t,x),\partial_x \f(t,x)\si(t),\cM[\f](t,x)\big)\ges g\big(t,x,\f(t,x),\partial_x \f(t,x)\si(t),\cM[u](t,x)\big).
\end{align*}
Thus, $\cL[\ti\f](t,x)\ges 0$, which yields that $\cL[ \f](t,x) \ges \cL[u](t,x)= 0$.
It means that $u(\cd,\cd)$ is a viscosity subsolution.
Moreover, we can show that $u(\cd,\cd)$ is also a viscosity supersolution.
\end{proof}

Under \ref{ass:H1}, the function $u(\cd,\cd)$ defined by \rf{u-Y} is merely continuous and does not
satisfy the regularity of classical solution.
In this case, we show that it is a viscosity solution under the additional monotonicity condition \ref{ass:H5}.

\begin{theorem}\label{viscosity-existence}
Under {\rm\ref{ass:H1}} and {\rm\ref{ass:H5}}, the function $u(\cd,\cd)$ defined by \rf{u-Y} is a viscosity solution to the MPE \rf{H-PDE-rw}--\rf{def-M-operator-rw}.
\end{theorem}

\begin{proof}
For any fixed $(t,x)\in[0,T]\times\dbR^n$, by the flow property of $Y^{t,x}(\cd)$, we have
\begin{align}
Y^{t,x}(t) &=Y^{t,x}(t+\d)+  \int_t^{t+\d} g(s, X^{t,x}(s),Y^{t,x}(s),Z^{t,x}(s),\cY^{t,x}(s)) ds- \int_t^{t+\d}Z^{t,x}(s) dW(s)\nn\\
&=Y^{t+\d,X^{t,x}(t+\d)}(t+\d)+  \int_t^{t+\d} g(s, X^{t,x}(s),Y^{t,x}(s),Z^{t,x}(s),\cY^{t,x}(s)) ds\nn\\
&\q- \int_t^{t+\d}Z^{t,x}(s) dW(s),\qq 0\les \d\les T-t.\nn
\end{align}
Thus
\begin{align}
u(t,x) &= u(t+\d,X^{t,x}(t+\d))+  \int_t^{t+\d} g(s, X^{t,x}(s),Y^{t,x}(s),Z^{t,x}(s),\cY^{t,x}(s)) ds\nn\\
\nn&\q- \int_t^{t+\d}Z^{t,x}(s) dW(s),\qq 0\les \d\les T-t.
\end{align}
Let $\f(\cd,\cd)$ be a smooth test function satisfying $\f(\cd,\cd) \ges u(\cd,\cd)$
and $\f(t, x) = u(t,x)$ for the fixed $(t,x)$.
Without loss of generality, we assume that $\f(\cd,\cd)$ has bounded derivatives.
Then
\begin{align}
\f(t,x) &\les \f(t+\d,X^{t,x}(t+\d))+  \int_t^{t+\d} g(s, X^{t,x}(s),Y^{t,x}(s),Z^{t,x}(s),\cY^{t,x}(s)) ds\nn\\
\nn&\q- \int_t^{t+\d}Z^{t,x}(s) dW(s),\qq 0\les \d\les T-t.
\end{align}
Notice that
\begin{align}
Y^{t,x}(r)
&= u(t+\d,X^{t,x}(t+\d))+  \int_r^{t+\d} g(s, X^{t,x}(s),Y^{t,x}(s),Z^{t,x}(s),\cM[u](s,X^{t,x}(s))) ds\nn\\
\nn &\q- \int_r^{t+\d}Z^{t,x}(s) dW(s),\qq t\les r\les t+\d.
\end{align}
Let $(\bar Y(\cd),\bar Z(\cd))$ be the unique solution to the BSDE:
\begin{align}
\bar Y^{t,x}(r) &= \f(t+\d,X^{t,x}(t+\d))+  \int_r^{t+\d} g(s, X^{t,x}(s),\bar Y^{t,x}(s),\bar Z^{t,x}(s),\cM[\f](s,X^{t,x}(s))) ds\nn\\
\nn&\q- \int_r^{t+\d}\bar Z^{t,x}(s) dW(s),\qq t\les r\les t+\d.
\end{align}
Note that
$$
\f(s,x)\ges u(s,x) \q\hbox{and}\q g(s, x,y,z,\cM[\f](s,x))\ges g(s, x,y,z,\cM[u](s,x)).
$$
Then by the comparison theorem of BSDEs, we have
$$
Y^{t,x}(s)\les\bar Y^{t,x}(s), \qq s\in[t,t+\d].
$$
Thus,
\begin{align}
\nn
\f(t,x)&=u(t,x)=Y^{t,x}(t)\les\bar Y^{t,x}(t) \\
\nonumber
&=\f(t+\d,X^{t,x}(t+\d))+  \int_t^{t+\d} g(s, X^{t,x}(s),\bar Y^{t,x}(s),\bar Z^{t,x}(s),\cM[\f](s,X^{t,x}(s))) ds\\
\nn
&\qq- \int_t^{t+\d}\bar Z^{t,x}(s) dW(s),\qq 0\les \d\les T-t.
\end{align}
Since $\f(\cd,\cd)$ is smooth, applying It\^o's formula to $s\mapsto\f(s,X^{t,x}(s))$ over $[t,t+\d]$, from the above we obtain
\begin{align}
\nn
&\int_t^{t+\d}\Big\{\partial_t\f(s,X^{t,x}(s))+{1\over 2}\tr\big [\partial_{xx}^2 \f(s,X^{t,x}(s))\si(s)\si(s)^\top\big]\nn\\
&\qq+g\big (s, X^{t,x}(s),\bar Y^{t,x}(s),\bar Z^{t,x}(s),\cM[\f](s,X^{t,x}(s))\big) \Big\}ds\nn\\
&- \int_t^{t+\d}\big[\bar Z^{t,x}(s)-\f_x(s, X^{t,x}(s))\si(s)\big]dW(s)\ges 0.\label{proof-vis-10}
\end{align}
Denote
\begin{align}
\nn
& \h Y^{t,x}(s)= \bar Y^{t,x}(s)-\f(s, X^{t,x}(s)),\qq  \h Z^{t,x}(s)= \bar Z^{t,x}(s)-\f_x(s, X^{t,x}(s))\si(s),\qq s\in[t,t+\d].\nn
\end{align}
Then
\begin{align}
\nn
\h Y^{t,x}(r)&=\int_r^{t+\d}\Big\{\partial_t\f(s,X^{t,x}(s))+{1\over 2}\tr[\partial_{xx}^2 \f(s,X^{t,x}(s))\si(s)\si(s)^\top]\nn\\
&\q+g\(s, X^{t,x}(s),\h Y^{t,x}(s)+\f(s, X^{t,x}(s)),\h Z^{t,x}(s)+\f_x(s, X^{t,x}(s))\si(s),\nn\\
&\qq\cM[\f](s,X^{t,x}(s))\) \Big\}ds- \int_r^{t+\d}\h Z^{t,x}(s)dW(s), \qq r\in[t,t+\d],\nn
\end{align}
and from \rf{proof-vis-10}, we have
\bel{proof-vis-13}
\h Y^{t,x}(t)\ges 0.
\ee
To calculate the limits of $(\h Y(\cd),\h Z(\cd))$ as $\d\to 0$, we introduce the following BSDE:
\begin{align}
\nn
\ti Y^{t,x}(r)&=\int_r^{t+\d}\Big\{\partial_t\f(s,x)+{1\over 2}\tr[\partial_{xx}^2 \f(s,x)\si(s)\si(s)^\top]+g\big(s, x,\ti Y^{t,x}(s)+\f(s, x),\nn\\
&\qq\ti Z^{t,x}(s)+\f_x(s, x)\si(s),\cM[\f](s,x)\big) \Big\}ds- \int_r^{t+\d}\ti Z^{t,x}(s)dW(s), \qq r\in[t,t+\d].\label{proof-vis-14}
\end{align}
Note that
$$
\dbE\[\sup_{s\in[t,t+\d]}|X^{t,x}(s)-x|^2\]\les \d.
$$
Then by the standard estimates of BSDEs, we have
$$
\dbE\[\sup_{s\in[t,t+\d]  }|\h Y ^{t,x}(s)-\ti Y^{t,x}(s)|^2\]\les K\dbE\[\(\int_t^{t+\d}|X^{t,x}(s)-x|ds\)^2\]\les K\d^3.
$$
In particular, we get
$$
|\h Y ^{t,x}(t)-\ti Y^{t,x}(t)|^2\les K\d^3.
$$
Combining the above with \rf{proof-vis-13}, we have
$$
{\ti Y^{t,x}(t)\over \d}\ges -K\delta^{1\over2}.
$$
%where ${\e(\d)\over \d}\to 0$ as $\d\to0$.
Obviously, the unique solution to \rf{proof-vis-14} is deterministic, so $\ti Z^{t,x}(\cd)\equiv 0$, and we have
\begin{align}
&{1\over \d}\int_t^{t+\d}\[\partial_t\f(s,x)+{1\over 2}\tr\big[\partial_{xx}^2 \f(s,x)\si(s)\si(s)^\top\big]\nn\\
&\q+g\big(s, x,\ti Y^{t,x}(s)+\f(s, x),\partial_x\f(s, x)\si(s),\cM[\f](s,x)\big) \]ds\ges -K\delta^{1\over2}.\label{proof-vis-18}
\end{align}
Applying the standard estimate of BSDEs, we have
$$
\dbE\[\sup_{s\in[t,t+\d]  }|\ti Y ^{t,x}(s)|^2\]\les K\d^2.
$$
Thus, letting $\d\to 0$ in \rf{proof-vis-18}, we have
\begin{align*}
&\partial_t\f(s,x)+{1\over 2}\tr\big[\partial_{xx}^2 \f(t,x)\si(t)\si(t)^\top\big]+g\big(t, x,\f(t, x),\partial_x\f(t, x)\si(t),\cM[\f](t,x)\big) \ges 0.\label{proof-vis-19}
\end{align*}
Thus the function $u(\cd,\cd)$ defined by \rf{u-Y} is a viscosity subsolution to the MPE \rf{H-PDE-rw}--\rf{def-M-operator-rw}.
Further, using the same arguments, we can show that $u(\cd,\cd)$ is also a viscosity supersolution to the MPE \rf{H-PDE-rw}--\rf{def-M-operator-rw}.
\end{proof}

Next we shall establish the comparison theorem for the viscosity solutions to MPE \rf{H-PDE-rw}--\rf{def-M-operator-rw}.
To begin with, we have a  partial comparison principle.

\begin{proposition}\label{thm-partial}
Suppose that {\rm\ref{ass:H1}} and {\rm\ref{ass:H5}} hold.
Let $u_1(\cd,\cd)$ and $u_2(\cd,\cd)$ be a viscosity subsoltion and supersolution to the MPE \rf{H-PDE-rw}--\rf{def-M-operator-rw}, respectively,
and denote $u(\cd,\cd):= u_1(\cd,\cd) - u_2(\cd,\cd)$.
Suppose
\bel{thm-partial-m1}
 \limsup_{|x|\to \infty} \sup_{0\les t\les T} {1\over 1+|x|^2}u(t,x) \les 0.
\ee
If $u_1(\cd,\cd)$ or $u_2(\cd,\cd)$ is smooth, then $u_1(\cd,\cd) \les u_2(\cd,\cd)$.
\end{proposition}

\begin{proof}
Without loss of generality, we assume that $u_2(\cd,\cd)$ is smooth.
Then  from \autoref{lem-consistency}, we have $\cL [u_2](t,x) \les 0$ for any $( t, x)\in[0,T]\times\dbR^n$.
Assume by contradiction that
$$
c_0 :=  \sup_{ 0\les t\les T,\, x\in \dbR^n} {1\over 1+|x|^2}u(t,x) >0.
$$
Set $\l := 2L+L^2+L\|M(\cd)\|_{L^\infty}+1$, where $L$ is the Lipschitz constant given in \ref{ass:H1}.
Then
$$
c_\l :=  \sup_{ 0\les t\les T,\, x\in \dbR^n}{1\over 1+|x|^2} e^{\l t}u(t,x) >0.
$$
By the boundary condition $u(T,x)=u_1(t,x)-u_2(t,x)\les 0$ and the assumption \rf{thm-partial-m1}, there exist $0\les t^*< T$ and $x^*\in\dbR^n$ such that
$${1\over 1+|x^*|^2}e^{\l t^*} u(t^*,x^*) = c_\l.$$
Define
$$
\f(t,x) := u_2(t,x)+c_\l(1+|x|^2) e^{-\l t}.
$$
Then $\f(\cd,\cd)\ges u_1(\cd,\cd)$ and  $\f(t^*,x^*)= u_1(t^*,x^*)$. Thus, $\f(\cd,\cd)$ is a test function of $u_1(\cd,\cd)$ at $(t^*, x^*)$, and then  $\cL [\f](t^*, x^*) \ges 0$.
Note that
\begin{align*}
&\partial_t \f (t,x)= \partial_t u_2(t,x) -\l c_\l(1+|x|^2) e^{-\l t},\qq \partial_x \f(t,x) = \partial_x u_2(t,x)+2c_\l xe^{-\l t},\\
& \partial_{xx} \f(t,x) = \partial_{xx} u_2(t,x)+2c_\l e^{-\l t}.
\end{align*}
Therefore,  at $( t^*, x^*)$,
\begin{align*}
0 &\les \cL[ \f] ( t^*, x^*)- \cL[u_2]( t^*, x^*) \\
&= - \l c_\l(1+|x^*|^2)  e^{-\l t^*}+c_\l e^{-\l t^*}\tr[\si(t^*) \si(t^*)^\top]\\
&\q+g\big(t^*, x^*,\f(t^*, x^*),\partial_x\f(t^*, x^*)\si(t^*),\cM[\f](t^*,x^*)\big)\\
&\q-g\big(t^*, x^*,u_2(t^*, x^*),\partial_x u_2(t^*, x^*)\si(t^*),\cM[u_2](t^*,x^*)\big)\\
&\les - \l c_\l(1+|x^*|^2)  e^{-\l t^*}+L^2 c_\l e^{-\l t^*}+L c_\l e^{-\l t^*}(1+|x^*|^2)+2L c_\l |x^*|e^{-\l t^*}\\
&\q+L \|M(\cd)\|_{L^\infty}c_\l(1+|x^*|^2){ e^{-\l t^*}\over\l} \\
&\les \[-\l +2L+L^2+{L \|M(\cd)\|_{L^\infty}\over\l}\] c_\l e^{-\l t^*}(1+|x^*|^2)\\
&\les -c_\l e^{-\l t^*}(1+|x^*|^2)<0,
\end{align*}
which is a contradiction. This completes the proof.
\end{proof}

\begin{remark}
Assumption \rf{thm-partial-m1} is a growth condition, and any continuous function $u(\cd,\cd)$ with at most linear growth in $x$ satisfies it.
\end{remark}

\ss
Since the regularity of $u(\cd,\cd)$ is only affected by the smoothness of the terminal condition $h(\cd)$
and the generator $g(\cd)$. We will establish the general comparison theorem by a modification method.

\begin{theorem} \label{thm:comparison-principle}
Suppose that {\rm\ref{ass:H1}} and {\rm\ref{ass:H5}} hold. Let $u_1(\cd,\cd)$ and $u_2(\cd,\cd)$ be any viscosity subsolution and supersolution to the MPE \rf{H-PDE-rw}--\rf{def-M-operator-rw}, respectively, both growing at most linearly in $x$. Then $u_1(\cd,\cd)\les u_2(\cd,\cd)$.
\end{theorem}

\begin{proof}
Let $u(\cd,\cd)$ be defined by \rf{u-Y},
it then follows from Theorem \ref{viscosity-existence} that $u(\cd,\cd)$ is a viscosity solution to the MPE \rf{H-PDE-rw}--\rf{def-M-operator-rw}.
We next prove the comparison principal under the Lipschitz continuous assumptions.
Define
\begin{align}
\nn h_\e(x)=h(x)-\e,\qq g_\e(s,x,y,z,y^\prime)=g(s,x,y,z,y^\prime)-\e,\\
\label{def-g-e} \forall (s,x,y,z,y^\prime)\in[0,T]\times\dbR^n\times\dbR\times\dbR\times\dbR.
\end{align}
By the mollification method, there exist smooth functions $\{h_{n,\e}(\cd),g_{n,\e}(\cd)\}_{n>0}$ such that
\begin{align}
&h_{n,\e}(x)\les h(x),\qq g_{n,\e}(s,x,y,z,y^\prime)\les g(s,x,y,z,y^\prime),\nn\\
\nonumber
&|g_{n,\e}(s,x_1, y_1,z_1,y^\prime_1)-g_{n,\e}(s,x_2, y_2,z_2,y^\prime_2)|\nn\\
\nn&\q\les L(|x_1-x_2|+|y_1-y_2|+|z_1-z_2|
+|y^\prime_1-y^\prime_2|),\\
\label{appro1}
&\qq\forall (s,x_i,y_i,z_i,y_i^\prime)\in[0,T]\times\dbR^n\times\dbR\times\dbR\times\dbR ,\q i=1,2,\q n>0,
\end{align}
and
\begin{align}
\nn\lim_{n\to\i }\big[ |h_{n,\e}(x)-h_{\e}(x)|+|g_{n,\e}(s,x,y,z,y^\prime)-g_\e(s,x,y,z,y^\prime)|\big]=0,\\
\label{appro2}  \forall (s,x,y,z,y^\prime)\in[0,T]\times\dbR^n\times\dbR\times\dbR\times\dbR.
\end{align}
Define
$$
u_{n,\e}(t,x)=Y_{n,\e}^{t,x}(t),
$$
where
$$
\left\{\begin{aligned}
&dY_{n,\e}^{t,x}(s)=-g_{n,\e}\big(s,X^{t,x}(s),Y_{n,\e}^{t,x}(s),Z_{n,\e}^{t,x}(s),\cM^{t,x}[Y_{n,\e}](s)\big)ds+Z_{n,\e}^{t,x}(s)dW(s),\qq s\in[t,T],\\
&Y_{n,\e}^{t,x}(T)=h_{n,\e}(X^{t,x}(T)),\end{aligned}\right.
$$
and for any fixed $s\in[t,T]$,
$$
\cM^{t,x}[Y_{n,\e}](s)=\int_s^T M(r-s) Y_{n,\e}^{r,X^{t,x}(s)}(r)dr,\qq s\in[t,T].
$$
Then by Theorem \ref{Thm-FK-classical}, $u_{n,\e}(\cd,\cd)$ is the unique classical solution to the following MPE:
$$
\left\{\begin{aligned}
&\partial_tu_{n,\e} (t,x)+{1\over 2}\tr\big[\partial_{xx}^2 u_{n,\e}(t,x)\si(t)\si(t)^\top\big]\\
&\q+g_{n,\e}\big(t,x,u_{n,\e}(t,x),\partial_x u_{n,\e}(t,x)\si(t),\cM[u_{n,\e}](t,x)\big)=0,\qq (t,x)\in[0,T]\times\dbR^n,\\
&u_{n,\e}(T,x)=h_{n,\e}(x),\q x\in\dbR^n,
\end{aligned}\right.
$$
where the memory operator $\cM$ is defined by
$$
\cM[u_{n,\e}](t,x)= \int_t^T M(s-t) u_{n,\e}(s,x)ds,\qq (t,x)\in[0,T]\times\dbR^n.
$$
Combining the above with \rf{appro1}, we get
$$
\cL [u_{n,\e}](t, x) \ges 0,\qq  u_{n,\e}(T, x)\les h(x),\qq \forall\,0 \les t\les T, \, x\in \dbR^n.
$$
Thus by \autoref{lem-consistency}, the smooth function $u_{n,\e}(\cd, \cd)$ is a viscosity subsolution to \rf{H-PDE-rw}--\rf{def-M-operator-rw}.
It follows from \autoref{thm-partial} that
\bel{u-n-2}
u_{n,\e}(t, x)\les u_2(t, x),\qq \forall\,0\les t\les  T,\, x\in \dbR^n.
\ee
Moreover, by \rf{def-g-e}, \rf{appro2} and the estimate of MBSDEs, we have
\bel{u-u-n}
\lim_{\e\to 0}\lim_{n\to\i}u_{n,\e}(t, x)= u(t, x),\qq \forall\,0\les t\les  T,\, x\in \dbR^n.
\ee
Combining \rf{u-n-2} with \rf{u-u-n}, we get
$$
u(t,x)\les u_2(t, x),\qq \forall\,0\les t\les  T,\, x\in \dbR^n.
$$
Similarly, we have
$$
u_1(t,x)\les u(t, x),\qq \forall\,0\les t\les  T,\, x\in \dbR^n.
$$
Thus,
$$
u_1(t,x)\les u_2(t, x),\qq \forall\,0\les t\les  T,\, x\in \dbR^n.
$$
This completes the proof.
\end{proof}

Combining Theorem \ref{viscosity-existence} and Theorem \ref{thm:comparison-principle}, we have the following result.

\begin{theorem}\label{viscosity-solution}
Under {\rm\ref{ass:H1}} and {\rm\ref{ass:H5}}, the function $u(\cd,\cd)$ defined by \rf{u-Y} is the unique viscosity solution
to the MPE \rf{H-PDE-rw}--\rf{def-M-operator-rw}.
\end{theorem}

\bs

\noindent
\textbf{Statement on no use of  AI:}
This manuscript was written entirely by the authors, and no AI-assisted tool  was used.

\linespread{1.0}

\end{document}